\documentclass[dvips,aop,preprint]{imsart}

\pdfoutput=1  %ADD FOR MATH ARXIV

\RequirePackage[OT1]{fontenc}

\RequirePackage{amsthm,amsmath, amssymb, amsbsy, latexsym, chicago, graphicx, color}
\RequirePackage[colorlinks,citecolor=blue,urlcolor=blue]{hyperref}
\RequirePackage{hypernat}

\arxiv{math.PR/1103.2317}

\startlocaldefs
\numberwithin{equation}{section}

\newcommand{\beq}{\begin{eqnarray}}
\newcommand{\eeq}{\end{eqnarray}}

\newcommand{\halmos}{\vspace{3mm} \hfill \mbox{$\Box$}}

\newcommand{\reals}{{\mathbb R}}

\newcommand{\integer}{{\mathbb Z}}
\newcommand{\pintegers}{{\integer_+}}
\newcommand{\integers}{{\mathbb N}}

\newcommand{\ra}{\rightarrow}
\newcommand{\ff}{\infty}
\newcommand {\Ed}{{\bf E}_{\mathfrak D}}
\newcommand{\Zt}{\tilde{Z}}
\newcommand{\Pf}{{\sc Proof.}\:\:}
\newtheorem{thm}{Theorem}[section]
\newtheorem{Prop}{Proposition}[section]
\newtheorem{Lm}{Lemma}[section]

\newtheorem{Rk}{Remark}[section]
\endlocaldefs

\def\Complex{\mathrm{\,\raise 0.33ex\hbox{\scriptsize\bf(}\!\!\!C}}
\def\Rational{\mathrm{\,\raise 0.33ex\hbox{\scriptsize\bf(}\!\!\!Q}}

\def\Bigbar#1{\mathrel{\left|\vphantom{#1}\right.\n@space}}

\begin{document}

\begin{frontmatter}
\title{Tail estimates for stochastic fixed point equations via  nonlinear renewal theory } %\protect\thanksref{T1}}
\runtitle{Tail estimates for homogeneous nonlinear recursions}
%\thankstext{T1}{Footnote to the title with the `thankstext' command.}

\begin{aug}
\author{\fnms{Jeffrey F.} \snm{Collamore}\thanksref{t1}\ead[label=e1]{collamore@math.ku.dk}} and
\author{\fnms{Anand N.} \snm{Vidyashankar}\thanksref{t2}\ead[label=e2]{avidyash@gmu.edu}}
%\ead[label=e3]{third@somewhere.com}
%\ead[label=u1,url]{http://www.foo.com}}

\thankstext{t1}{Corresponding author.  Research supported in part by Danish Research Council (SNF) Grant ``Point Process
Modelling and Statistical Inference," No. 09-092331.}
\thankstext{t2}{Research Supported by grants from NDC Health Corporation and 
NSF DMS 000-03-07057.}
%\thankstext{t3}{Second supporter of the project}
\runauthor{J.F.\ Collamore and A.N.\ Vidyashankar}

\affiliation{University of Copenhagen and George Mason University}

\address{Department of Mathematical Sciences\\
University of Copenhagen\\
Universitetsparken 5\\
DK-2100 Copenhagen {\O}\\
Denmark\\
\printead{e1}\\
\phantom{E-mail:\ }\printead*{}}

\address{Department of Statistics\\
George Mason University\\
4400 University Drive,
MS 4A7\\
Fairfax, VA  22030\\
U.S.A.\\
\printead{e2}\\
\printead{}}
\end{aug}

\begin{abstract}
This paper presents precise large deviation estimates
for solutions to stochastic fixed point equations
of the type $V \stackrel{\cal D}{=} f(V)$, where $f(v) = Av + g(v)$
for a random function $g(v) = o(v)$ a.s.\ as $v \to \infty$.
Specifically, we provide an explicit characterization of the
pair $(C,\xi)$ in the tail estimate
${\bf P} \left( V > u \right) \sim C u^{-\xi}$ as $u \to \infty$,
and also present a Lundberg-type upper bound of the form ${\bf P} \left( V > u \right) \le \bar{C}(u) u^{-\xi}$.  To this end, we introduce a novel dual
change of measure on a random time interval and analyze the path
properties, using nonlinear renewal theory, of the Markov chain
resulting from the forward iteration of the given stochastic fixed
point equation.  In the process, we establish several new results
in the realm of nonlinear renewal theory
for these processes.  As a consequence of our techniques,
we also establish a new characterization of the extremal index.
Finally, we provide some extensions of our methods to Markov-driven
sequences.
\end{abstract}

\begin{keyword}[class=AMS]
\kwd[Primary ]{60H25}
\kwd[; secondary ]{60K05}
\kwd{60F10}
\kwd{60J10}
\kwd{60G70}
\kwd{60K25}
\kwd{60K35}.
\end{keyword}

\begin{keyword}
Random recurrence equations, Letac's principle, nonlinear renewal theory, slowly changing functions, Harris recurrent
Markov chains, geometric ergodicity, large deviations,
Cram\'{e}r-Lundberg theory with stochastic investments, GARCH processes, extremal index.
\end{keyword}

\end{frontmatter}

\section{Introduction}
\subsection{Motivation and statement of the problem}
Stochastic fixed point equations (SFPE) arise in several areas of contemporary science and have been the focus of much study in applied probability, finance, analysis of algorithms, page ranking in personalized web search,
risk theory and actuarial mathematics. A general stochastic fixed point equation has the form
\beq{\label{eq0.1}}
V \stackrel{\cal{D}}=f(V),
\eeq
where $f$ is a random function satisfying certain regularity
conditions and is independent of $V$. When $f(v)= Av +B$, where $(A, B)$ is independent of $V$ and ${\bf E} \left[ \log A \right] <0$, this problem has long history starting with the works of \shortciteN{FS72}, \shortciteN{HK73},  \shortciteN{WV79}, and \shortciteN{GL86}.

Tail estimates for solutions of general SFPEs have been developed by \shortciteN{CG91} using implicit renewal theory, extending the fundamental
work of \shortciteN{HK73}. 
Under appropriate regularity conditions, Goldie established that 
\beq{\label{eq0.2}}
{\bf P} \left( V > u \right) \sim C u^{-\xi} \quad \mbox{\rm as} \quad u \ra \ff,
\eeq
where here and in the following, $g(u) \sim h(u)$ means $\lim_{u \to
  \infty} g(u)/h(u) = 1$. The constant $C$ appearing in (\ref{eq0.2})
is defined in terms of the tails of $V$ and $f(V)$, rendering it
unsuitable for numerical purposes and for statistical inference
concerning various parameters of scientific interest
that involve the tails of $V$. Indeed, except in some
special cases, the formulae for $C$ presented in \shortciteN{CG91} do
not simplify to useful expressions.  This issue is somewhat folklore
and has been observed by several researchers
including \shortciteN{BYMP98}, 
\shortciteN{DS01}, and \shortciteN{NECSOZ09}.
Yakin and Pollak's work was motivated by likelihood ratio testing and
change point problems in statistics, while 
the paper of Enriquez et al. was motivated by probabilistic considerations.

The primary objective of this paper is to present a general
alternative probabilistic approach to deriving the above tail
estimates, yielding a characterization of this constant which---beyond
its theoretical interest---is
also amenable to statistical
inference, and Monte Carlo estimation in particular.
In the process, we draw an interesting
connection between the constant $C$ in \eqref{eq0.2} and the backward
iterates of an SFPE described in \shortciteN{GL86}.  Specifically, we show that
in many cases, this constant may be obtained via iterations of a
corresponding perpetuity sequence, which has a comparatively simple form and is computable.
As this representation is {\em explicit}, it also resolves questions
about the positivity of the constant which have been often raised in
the literature and addressed via {\it ad hoc} methods.  In addition,
based on a slight variant of our method,
we develop a sharp upper bound for the the tails of 
$V$, which has a form analogous to the classical Lundberg inequality from
actuarial mathematics.  

The key idea to our approach is the observation that the  distribution of the solution to the SFPE is frequently the stationary limit of a positively recurrent Markov chain 
and, hence, tail probabilities can be studied via excursions within a regeneration cycle of this chain.  After a large deviation change of measure argument, these excursions can be viewed as perturbations of a multiplicative random walk. While the multiplicative random walk itself can be handled via classical renewal-theoretic techniques, analysis of perturbations requires the development of new methods in nonlinear renewal theory. 

The methods of this paper are quite general and apply to a wide class of SFPEs, including the following SFPEs which will be analyzed 
in some detail below:
\begin{eqnarray}
{\label{eq0.3A}}V  &\stackrel{\cal D}=& \left( AV + B \right)^+,\\%~ V ~ \mbox{independent of }~ (A, B), \\
{\label{eq0.3C}}V &\stackrel{\cal D}=& A \max (D, V) + B,\\ 
{\label{eq0.3D}}V &\stackrel{\cal D}=& AV + BV^{1/2} +D.
\end{eqnarray}
In the above,  $V$ is a real-valued random variable and independent of the other random variables on the right-hand side of these equations, while the random variables in the pair $(A,B)$ or triplet $(A,B,D)$ are allowed to be dependent on one another.
Under the assumption that ${\bf E} \left[ \log A \right] < 0$, a unifying feature of these equations is that they may all be written in the rough form $V \stackrel{\cal D}{=} AV + $``error," where the error term becomes negligible in the asymptotic limit for large $V$.

It is also of considerable interest to study the case where ${\bf E}
\left[ \log A \right] \ge 0$. Borrowing the terminology from
autoregressive processes, one may regard the case where ${\bf E}
\left[ \log A \right] > 0$ as the  {\it explosive case}  while the case
where ${\bf E} \left[ \log A \right] < 0$  as the {\it stationary case}.
In a sequel to this work, we establish central limit theorems and
rates of convergence results under appropriate scalings for the
explosive case. Indeed, we establish that under some specific
conditions, $(V_n/\rho^n) \to W \quad \mbox{\rm
  as} \quad n \to \infty$ with probability 1 (w.p.1) for a certain constant $\rho$. This limit
random variable $W$ will satisfy an SFPE, and the results of this
article can be invoked to provide sharp tail estimates for
$W$. Finally, the case where ${\bf E} \left[ \log A \right]=0$,
referred to as the critical case in the literature, has been studied
in detail in \shortciteN{BMBPEL97} and \shortciteN{BD07} under fairly
stringent conditions on the random variables $B$ and $D$. However,
methods of the current paper can be adopted to reduce some of these
conditions. Hence, combining the results of this paper with known, but
modified, results for the critical case and our results for the
explosive case, one can obtain a {\em complete}\, description of the
tail behavior for a wide class of random variables satisfying an SFPE encompassing \eqref{eq0.3A}-\eqref{eq0.3D}.

We now turn to a more detailed description of specific examples that originally motivated this work.

%
%
%  RUIN WITH INVESTMENTS
%
%

\subsection{The ruin problem with stochastic investments}
A fundamental problem in actuarial mathematics---first introduced by \shortciteN{FL03} and \shortciteN{HC30}---considers the probability
of eventual ruin for an insurance company.  In this problem, the company is assumed to receive premiums at a constant rate $c$
and incur claim losses which arise according to a Poisson rate.  Thus, the capital growth of the insurance company is modeled as
\begin{equation} \label{cl1}
X_t = c t - \sum_{i=1}^{N_t} \zeta_i,
\end{equation}
where $\{N_t \}$ is a Poisson process and the loss sequence $\{
\zeta_i \}$ is assumed to be independent and identically distributed (i.i.d.). The ruin problem is concerned with the probability, starting from a positive initial capital of $u$, that the company ever incurs ruin.  Since $X_0=0$, 
this is equivalent to
\begin{equation}
\Psi^\ast(u) := {\bf P} \left( X_t < -u, \:\: \mbox{\rm for some} \:\: t \ge 0 \right).
\end{equation}
Cram\'{e}r's classical result states 
\begin{equation} \label{cl2}
\Psi^\ast(u) \sim C^\ast e^{-\xi u} \quad \mbox{\rm as} \quad u \to \infty,
\end{equation}
for certain positive constants $C^\ast$ and $\xi$, 
while the well-known Lundberg inequality states
\begin{equation} \label{cl3}
\Psi^\ast(u) \le e^{-\xi u} \quad \mbox{\rm for all} \:\: u \ge 0.
\end{equation}
The constant $\xi$ is called the Lundberg exponent and the constant $C^{\ast}$ is explicitly characterized in (\ref{sr203}) below.

A natural criticism of Cram\'er's model is that the company's assets are assumed to be invested in a riskless bond yielding zero interest, while in practice these assets would typically be invested in a portfolio of stocks and bonds earning positive interest. One method to address this shortcoming is to introduce a stochastic investment process.
To this end,  let $ R_n $ denote the returns on the investments in the $n$th discrete time interval, and let  $ L_n $ denote the losses from the insurance business during the same time interval; that is,
\[
L_n = X_{n} - X_{n-1},
\]
where $\{ X_t \}_{t \ge 0}$ is as given in \eqref{cl1}.  The total capital $Y_n$ of the company at time $n$ then satisfies the recursion
\begin{equation} \label{intro2}
Y_n := R_n Y_{n-1} - L_n,  \quad n = 1,2,\ldots, \quad Y_0 = u.
\end{equation}
 The ruin probability is now given by
$\Psi(u) := {\bf P} \left( Y_n < 0, \mbox{ some }n \right)$,
and by an elementary argument it can be seen that 
\begin{equation} \label{intro3}
\Psi(u) = {\bf P} \left( {\cal L}_n > u, \mbox{ some } n \right),
\end{equation}
where
\begin{equation} \label{intro3a}
{\cal  L}_n := R_1^{-1} L_1 + R_1^{-1}R_2^{-1} L_2 + \cdots + (R_1 \cdots R_{n})^{-1} L_n.
\end{equation}
Using that
$R_2^{-1}L_2 +\cdots +(R_2 \cdots R_{n+1})^{-1}L_{n+1} \stackrel{\cal{D}}={\cal L}_n,$
it follows from \eqref{intro3a} that
\begin{equation} \label{intro4a}
{\cal L}_n \stackrel{\cal{D}}= B + A {\cal  L}_{n-1} \quad \mbox{\rm
  where} \quad A= \frac{1}{R_1} \:\: \mbox{\rm and} \:\: B= \frac{L_1}{R_1}.
\end{equation}
Hence, setting $\left( \sup_{n \in {\mathbb Z}} {\cal L}_n \right)
 \vee 0 = {\cal L}$, we obtain the SFPE
\begin{equation} \label{cl10}
{\cal L} \stackrel{\cal D}{=} \left( A {\cal L} + B \right)^+,
\end{equation}
which is  \eqref{eq0.3A}.
Appealing to \shortciteN{CG91}'s implicit renewal theorem yields
\begin{equation}\label{intro5}
\Psi(u) \sim C u^{-\xi} \quad \mbox{as} \quad u \to \infty,
\end{equation}
obtaining a polynomial decay for ruin, which should be contrasted with the exponential decay in Cram\'{e}r's original estimate.

Variants on Goldie's result have been extensively studied in the actuarial literature; see, {\it e.g.}, \shortciteN{JP08} and the references therein. 
For an extension to Harris recurrent Markov chains
in general state space, see \shortciteN{JC09}.  In spite of this extensive literature, we are not aware of a canonical method for
characterizing the constant $C$ in this and a variety of similar applications, and our paper fills this void.
In practice,  statistical estimates for $C$ are of equal importance
 as those for $\xi$, and  \eqref{intro5} is not very useful without a method for estimating these constants.  Furthermore, as a by-product of our approach, we also obtain an analog of the classical Lundberg inequality for this problem, namely,
\begin{equation} \label{cl15}
\Psi(u) \le \bar{C}(u)  u^{-\xi} \quad \mbox{\rm for all}\:\: u \ge 0,
\end{equation}
where $\bar{C}(u)$ converges to a positive constant as $u \to \infty$.

There is an interesting connection of this problem to the extremes of the well-known ARCH(1) and GARCH(1,1) financial models, which satisfy
a related SFPE, namely $V \stackrel{\cal D}{=} AV + B$.
We will pursue this connection in more detail in Section 3, where we also apply the methods of this paper to yield a new, explicit characterization
for the much-studied extremal index.

%, as illustrated in Figure 1.
%
%
%  FIGURE 1
%
%
%\begin{figure}[f1]
%\begin{center}
%\subfigure[]{
%\includegraphics[scale=0.4]{compare_norms}}
%\subfigure[]{
%\includegraphics[scale=0.4]{CSlargep100rhoff}}
%\caption{Kernel density estimates of the norms of statistics and the Histogram of p-values under the sup norm.}
%\caption{\label{fig:result} }
%\end{center}
%\end{figure}

%\begin{figure}[f1]
%\begin{center}
%\includegraphics{Figure1.pdf}
%\caption{A sample path of the classical Cram\'{e}r-Lundberg process with initial capital $u=5$.  Ruin occurs when the process
%$\{ X_t \}$ falls below zero.}
%\label{Fig1}
%\end{center}
%\end{figure}

%One common reason that is usually cited is that the risk process $\{ X_t \}$ may not be homogeneous in time;
%addressing of this issue has  lead to numerous extensions of the basic model,  as summarized in \shortciteN{SA00}.

%A sample path for this cumulative loss process is shown in Figure 2, and should be compared with the classical Cram\'{e}r-Lundberg
%process in Figure 1.
%
%
%  FIGURE 2
%
%

%\begin{figure}[f2]
%\begin{center}
%\includegraphics[scale=0.2]{Figure2.pdf}
%\caption{A sample path of the cumulative loss process.  Ruin occurs when this negative-drift process reaches a positive barrier at $u$, where $u$ is the initial
%capital.  Note that the process has higher volatility than the classical Cram\'{e}r-Lundberg process.
%Moreover, as seen from the comparatively long time scale, the jumps sizes typically become smaller with increasing $t$.}
%\label{Fig2}
%\end{center}
%\end{figure}

%
%
%  LETAC'S MODEL E
%
%

\subsection{Related recursions and Letac's Model E}
The models described in the previous section are special cases of a more general model introduced in \shortciteN{GL86}.
Specifically, Letac's Model E is obtained via the forward recursion
\begin{equation} \label{eq0.12}
V_n = A_n\max(D_n, V_{n-1}) + B_n, \quad n=1,2,\ldots \quad \mbox{\rm and}\quad V_0= v,
\end{equation}
where the {\em driving sequence}
 $\{ (A_n,B_n,D_n):  n=1,2,\ldots \}$ is i.i.d., and we will assume in this discussion and throughout the paper that $\{ A_n \}$ is positive-valued and ${\bf E} \left[ \log A \right] < 0.$
Iterating the above equation leads to the  expression
\begin{equation}{\label{eq0.13}}
V_n =\max \left(\ \sum_{i=0}^n B_i \prod_{j=i+1}^n A_j, \bigvee_{k=1}^n \left[  \sum_{i=k}^n B_i \prod_{j=i+1}^n A_j + D_k \prod_{j=k}^n A_j
\right] \right)
\end{equation}
with $B_0 \equiv v$.
The sequence $V_n$ is a Markov chain, and one can show that under moment and regularity conditions, $V: = \lim_{n \to \infty} V_n$ exists almost
surely (a.s.).  
A more compact expression is obtained in the form of an SFPE upon observing that \eqref{eq0.12} implies
\begin{equation}{\label{eq0.14}}
V \stackrel{\cal D}{=} A \max(D,V) + B,
\end{equation}
which is \eqref{eq0.3C}.
As we will discuss in more detail in the next section, there is also a
corresponding backward iteration.  Setting $B_i^\ast = B_i$ for
$i=1,\ldots,n$ and $B_{n+1}^\ast = v$, then
iteration of the backward equation yields that
\begin{equation} \label{eq0.15}
Z_n = \max \left( \sum_{i=1}^{n+1} B_i^\ast \prod_{j=1}^{i-1} A_j,
\bigvee_{k=1}^n \left[ \sum_{i=1}^k B_i \prod_{j=1}^{i-1}A_j + D_k \prod_{j=1}^k A_j
\right] \right).
\end{equation}
In fact, it can be shown that both the forward and the backward recursions satisfy the same
SFPE, namely \eqref{eq0.14}.  

Letac's Model E is quite general, encompassing a wide variety of recursions appearing in the literature.  For example, setting 
$D=-B/A$ in \eqref{eq0.14} yields $V=(AV+B)^+$, which leads us back to
the SFPE arising in our discussion of the ruin problem with
investments. Similarly, the recursions arising in the ARCH(1) and
GARCH(1,1) processes are easily seen to be special cases of
\eqref{eq0.14} with $D=0$ and $B$ strictly positive.
Furthermore, we can also obtain the so-called
perpetuity sequences arising in life insurance mathematics.  Let
$ B_n $ be nonnegative, $V_0=0$, and $D_n = 0$ for all $n$; then the backward recursion \eqref{eq0.15}
simplifies to
\begin{equation} \label{eq0.16}
Z_n = B_1 + A_1 B_2 + A_1 A_2 B_3 + \cdots + (A_1 \cdots A_{n-1}) B_n, \quad n=1,2,\ldots.
\end{equation}
In a life insurance context, $ B_n $ would represent future liabilities to the company, and the quantity $\lim_{n \to \infty} Z_n$ would then describe the discounted future liabilities faced by the company. The random variables $ A_n $ are the ``discount factors" which, as in the ruin problem with investments, may be viewed as the inverse returns of an appropriate financial process.

Letac's Model E also includes reflected random walk, which arises
in classical queuing theory with the study of the steady-state waiting
times for a GI/G/1 queue.  Indeed, setting $B=0$ and $D=1$ in \eqref{eq0.14}, we obtain an SFPE corresponding to the forward equation
\begin{equation}
V_n = \max \Big(1, \bigvee_{j=1}^n \prod_{i=j}^n A_i \Big),
\end{equation}
which is seen to correspond with reflection at $ 1$ of the multiplicative random walk $\prod_{i=1}^n A_i$ or, equivalently,
to reflection at $ 0 $ of the random walk $S_n := \sum_{i=1}^n \log A_i$.
Using classical duality,  the exceedances of $\{ V_n \}$ can be related
to maxima of the unreflected random walk.
Similar arguments lead to further connections with classical risk models, autoregressive processes with random coefficients, 
and related problems.\\[-.25cm]

The rest of the paper is organized as follows.  Section 2 contains notation, assumptions, background on Letac's principle and connections with
Markov chain theory, and the main theorems of the paper. Section 3 is devoted examples and consequences of our main theorems. In Section 4 we describe our
results in nonlinear renewal theory. These results serve as the essential technical tools that are needed for the proofs of the main results of this paper and seem to be of independent interest.  
Section 5 is devoted to proofs of the main theorems stated in Sections 3, while Section 6 contains the proofs of results on nonlinear renewal theory.  An extension of Letac's Model E to nonlinear models is described in Section 7,
while Section 8 is devoted to extensions for the Markov case, and Section 9 
to some concluding remarks.

%
%
%  STATEMENT OF RESULTS
%
%
\section{Statement of results}
\subsection{Letac's principle and background from Markov chain theory}
In the previous section, we presented forward and backward equations associated with the SFPE \eqref{eq0.14}.  In particular,
the forward equation is generated by converting the SFPE \eqref{eq0.14} into the recursion \eqref{eq0.12}, while the backward equation
proceeds in similar fashion, but using an iteration backward in time. 

More precisely, suppose that we begin with the general SFPE
\begin{equation} \label{letac1a}
V \stackrel{\cal D}{=} f(V),
\end{equation}
and assume that this equation can be written as a deterministic functional equation of the form
\begin{equation} \label{letac1}
V \stackrel{\cal D}{=} F(V,Y),
\end{equation}
where $F:  \reals \times \reals^d \to \reals$ is assumed to be deterministic, measurable, and  continuous in its first component.  For convenience of notation,
we will often write $F_y(v) = F(v,y)$, where $y \in \reals^d$.   For example, the SFPE \eqref{eq0.14} may be written
as $V \stackrel{\cal D}{=} F_Y(V)$, where $F(v,Y) = A \max \left( D,v \right) + B$ and $Y=(A,B,D)$.

Let $v $ be an element of the range of $F$, and let $\{ Y_n \}$ be an i.i.d.\ sequence of random variables such that $Y_i \stackrel{\cal D}{=} Y$ for all $i$.
Then the forward equation generated by the SFPE \eqref{letac1} is defined by
\begin{equation} \label{letac2}
V_n (v) = F_{Y_n} \circ F_{Y_{n-1}} \circ \cdots \circ F_{Y_1}(v), \quad n=1,2,\ldots, \quad V_0 = v,
\end{equation}
while the backward equation generated by this SFPE is defined by
\begin{equation}
Z_n(v) = F_{Y_1} \circ F_{Y_2} \circ \cdots \circ F_{Y_n} (v), \quad n=1,2,\ldots, \quad Z_0 = v.
\end{equation}

Note that the backward recursion need not be Markovian, but for every
$v$ and every integer $n$,
 $V_n(v)$ and $Z_n(v)$ are identically distributed. Letac's (1986)
 principle states that  when the backward recursion converges to $Z$
 a.s. and is independent of its starting value $v$,
 then the stationary distribution
 of $V_n$ is
the same as the 
distribution of $Z$. This is precisely described in the following lemma.

%
%
%  LETAC'S PRINCIPLE
%
%

\begin{Lm} \label{LetacsLM}
Let  $F:  \reals \times \reals^d \to \reals$ be deterministic, measurable, and  continuous in its first component; that is,
$F_Y$ is continuous, where $Y:  \reals^d \to \reals$ for some 
$d \in \pintegers$.  Suppose that $\lim_{n \to \infty} Z_n(v) := Z$ exists
a.s.\ and is
independent of its initial state $v$.   Then $Z$ is the unique
solution to the SFPE \eqref{letac1}. Furthermore, $V_n(v)
\Rightarrow V$ as $ n \ra \ff$, where $V$ is independent of $v$, and
the law of $V$ is same as the law of $Z$.
\end{Lm}

Starting from the general SFPE \eqref{letac1}, our basic approach will be to generate the forward recursive sequence \eqref{letac2}
and study the evolution of this process.
In the next section, we will impose additional assumptions that will ensure that the forward recursive sequence viewed as a Markov chain possesses a  regenerative structure, which we now describe.

For the remainder of this section, let $\{ V_n \}$ denote a
Markov chain on a general state space $({\mathbb S},{\cal S})$,
and assume that this chain is aperiodic, countably generated, and irreducible with respect to its maximal irreducibility measure, which we will denote by $\varphi$.  Let $P$ denote the transition kernel of $\{ V_n \}$. In Section 2.2 onward, we will take $\{ V_n \}$ to be the Markov chain generated by a forward recursion, but for the remainder
of this section we will assume it to be an arbitrary chain.

Recall that if  $\{ V_n \}$ is aperiodic and $\varphi$-irreducible, then it
satisfies a {\em minorization}, namely
\begin{equation}
\delta {\bf 1}_{\cal C} (x) \nu(E) \le
P^k(x,E),
\quad\forall x \in {\mathbb S},\: E \in {\cal S}, \quad k \in \pintegers , \tag{${\cal M}$}
\end{equation}
for some set ${\cal C}$ with $\varphi({\cal C})>0$, some constant $\delta > 0$, and some probability measure
$\nu$ on $({\cal S}, {\mathbb S})$.   The set ${\cal C}$ is called the ``small set" or ``${\cal C}$-set."
It is no loss of generality to assume that this ${\cal C}$-set is bounded,
and we will assume throughout the article that this is the case.

An important consequence of $({\cal M})$, first established in \shortciteN{KAPN78} and \shortciteN{EN78},
is that it provides a regeneration structure for the chain.  More precisely, these authors established the following.

\begin{Lm} \label{regnLM}
Assume that $({\cal M})$ holds with $k=1$.
Then there exists a sequence of random times, $0\le K_0< K_1 < \cdots$, such that:\\
\hspace*{.23in}{\rm (i)} $K_0, K_1 -K_0, K_{2} - K_{1},\ldots$
are finite a.s.\ and mutually independent{\rm ;}\\
\hspace*{.19in}{\rm (ii)} the sequence $\{ K_{i} - K_{i-1} :
i=1,2,\ldots \}$ is i.i.d.{\rm ;}\\
\hspace*{.18in}{\rm (iii)} the random blocks $\left\{ V_{K_{i-1}},\ldots,V_{K_{i}-1} \right\}$ are independent,
$i=0,1,\ldots${\rm ;}\\
\hspace*{.18in}{\rm (iv)} ${\bf P}\left(V_{K_i} \in E | {\mathfrak F}_{K_i-1}\right) = \nu (E)$, for all $E \in
{\cal S}$.
\end{Lm}
\noindent
Here and in the following, ${\mathfrak F}_n$ denotes the $\sigma-$field generated by $V_0, V_1, \ldots, V_n$.
Let $\tau_i := K_i - K_{i-1}$ denote
the $i$th inter-regeneration time, and let
$\tau$ denote a typical regeneration time, {\it i.e.}, $\tau \stackrel{\cal D}{=} \tau_1$.

\begin{Rk}{\rm
Regeneration can be related to the return times of $\{ V_n \}$ to the ${\cal C}$-set in $({\cal M})$ by introducing an augmented
chain $\{ (V_n,\eta_n) \}$,  where $\{ \eta_n \}$ is an i.i.d.\ sequence of Bernoulli random variables, independent of $\{ V_n \}$
with ${\bf P} \left( \eta_n = 1 \right) = \delta$.  Then based on the proof of Lemma 2.2, we may identify $K_{i}-1$ as the $(i+1)$th
return time of the chain $\{ (V_n,\eta_n) \}$ to the set ${\cal C} \times \{ 1 \}$.  At the subsequent time, $V_{K_i}$ then has the distribution $\nu$
given in $({\cal M})$, independent of the past history of the Markov chain.  Equivalently, the chain $\{ V_n \}$ regenerates with probability
$\delta$ upon each return to the ${\cal C}$-set.  In particular, in the special case that $\delta = 1$, in which case the chain $\{ V_n \}$ is
said to have an atom, the process $\{ V_n \}$ regenerates upon {\em every} return to the ${\cal C}$-set.
}
\end{Rk}

%
%
%  REMARK 2.2:  Removed!
%
%
%
%\begin{Rk}{\rm
%In $({\cal M})$ there is no loss of generality in assuming that the
%set $C$ and the support of $\nu$ are bounded, 
%and in th% sequel we will always assume
% this to be the case.}
%\end{Rk}

Note that if the Markov chain is also recurrent, then $\tau < \infty$ a.s.
However, in this article, we will also consider transient chains which
may {\em never} regenerate.  In that case, we will set $\tau = \infty$.

As consequence of the previous lemma,
one obtains a representation formula (\shortciteN{EN84}, p.\ 75) relating the stationary limit distribution of $V$
%\[
%V := \lim_{n \to \infty} V_n
%\]
to the behavior of the chain over a regeneration cycle.  Set
\begin{equation} \label{defNu}
N_u = \sum_{n=0}^{\tau -1} {\bf 1}_{(u,\infty)}(V_n),
\end{equation}
where $V_0 \sim \nu$.  Thus, using Lemma 2.2 (iv), we see
that $N_u$ describes the number of visits of $\{ V_n \}$ to the set $(u,\infty)$ over a typical regeneration cycle. The representation formula alluded to above is given in the following lemma.

\begin{Lm}  Assume that $({\cal M})$ holds with $k=1$.  Then for any $u \in \reals$,
\begin{equation} \label{MC-1}
{\bf P} \left( V > u \right) = \frac{{\bf E} \left[ N_u \right]}{{\bf E} \left[ \tau \right]}.
\end{equation}
\end{Lm}

%For further background on Markov chain theory, see, {\it e.g.}, \shortciteN{EN84} or \shortciteN{SMRT93}.

\subsection{Results in the setting of Letac's Model E}
We now turn to the core results of the paper.   In this discussion, we
will focus primarily on Letac's Model E, namely,
\begin{equation} \label{sr200}
V \stackrel{\cal D}{=}  F_Y(V), \quad \mbox{\rm where}\quad   F_Y(v) =  A \max(D,v) + B
\end{equation}
for $Y = (A,B,D) \in (0,\infty) \times \reals^2$ (although we will
later generalize the method).
As described previously,
this model is sufficiently general to incorporate the main
applications described in the introduction.  Let $\{V_n\}$ and $\{ Z_n
\}$ denote the forward and backward recursions generated by this SFPE,
respectively.
In addition to these sequences, we also introduce two related backward
recursive sequences, which
will play an essential role in the sequel.

%
%
%  ASSOCIATED PERPETUITY SEQUENCE
%
%

{\it The associated perpetuity sequence}{\rm :}  Let $\{ (A_n,B_n,D_n) \}_{n
  \in \pintegers}$
be the i.i.d.\ driving sequence which generates the forward recursion $\{
V_n \}$, and let $A_0 = 1$ and $B_0 \sim \nu$ for $\nu$ given as in
Lemma 2.2 (iv), where $B_0$ is independent of the driving sequence $\{
(A_n,B_n,D_n) \}_{n \in \pintegers}$.
Now consider the backward recursion
\begin{equation} \label{def-aps1}
Z_n^{(p)} = F^{(p)}_{Y_0} \circ \cdots \circ F^{(p)}_{Y_n}(0), \quad n=0,1,\ldots,
\end{equation}
where
\begin{equation} \label{def-aps2}
F^{(p)}_Y(v) := \frac{v}{A} + \frac{B}{A}.
\end{equation}
By an elementary inductive
argument, it follows that
\begin{equation} \label{aps3}
Z_n^{(p)} = \sum_{i=0}^n \frac{B_i}{A_0 \cdots A_i}, \quad n=0,1,\ldots.
\end{equation}
The sequence $\{ Z_n^{(p)} \}$ will be called the {\em perpetuity sequence
associated to} $\{ V_n \}$ and is the backward recursion
generated by the SFPE \eqref{def-aps2}.

%
%
%  CONJUGATE SEQUENCE
%
%

{\it The conjugate sequence}{\rm :}  
Let $D_0^\ast := -B_0$ and $D_i^\ast := -A_i D_i - B_i$ for
$i=1,2,\ldots,$
and consider the backward recursion
\begin{equation} \label{def-acs1}
Z_n^{(c)} = F^{(c)}_{Y_0} \circ \cdots \circ F^{(c)}_{Y_n}(0), \quad n=0,1,\ldots,
\end{equation}
where
\begin{equation} \label{def-acs2}
F^{(c)}_Y(v)  := \frac{1}{A} \min \left( D^\ast,v \right) +
\frac{B}{A}.
\end{equation}
It follows by induction that
\begin{equation} \label{def-acs4}
Z_n^{(c)} = \min \left(Z_n^{(p)}, 0 , 
  \bigwedge_{k=1}^n \sum_{i=0}^{k-1} \frac{B_i}{A_0 \cdots A_i}
- \frac{D_k}{A_0 \cdots A_{k-1}}
 \right).
\end{equation}
The sequence $\{ Z_n^{(c)} \}$ will be called the {\em conjugate sequence
to} $\{ V_n \}$ and is the backward recursion
generated by the SFPE \eqref{def-acs2}.

{\it Further notations}{\rm :}  Let $A$ denote the multiplicative factor appearing in \eqref{sr200}, and define
\[
\lambda (\alpha) = \log {\bf E} \left[ A^\alpha \right] \quad \mbox{and} \quad
\Lambda(\alpha) = \log \lambda(\alpha), \quad \forall \alpha.
\]

Let $\mu$ denote the distribution of $Y:= (\log A,B,D)$, and let $\mu_\alpha$ denote the $\alpha$-shifted
distribution with respect to the first variable; that is,
\begin{equation} \label{meas-shift}
\mu_\alpha(E) := \int_E e^{\alpha x}Êd\mu(x,y,z), \quad \forall E \in {\cal B}(\reals^3), \quad \forall \alpha \in \reals,
\end{equation}
where here and in the following, ${\cal B}(\reals^d)$ denotes the Borel sets of $\reals^d$ for any $d \in \pintegers$.
Let ${\bf E}_\alpha\left[ \cdot \right]$ denote expectation with respect to this $\alpha$-shifted measure. 
For any random variable $X$, let
 ${\mathfrak L}(X)$ denote the probability law of $X$.  Also, write $X \sim {\mathfrak L}(X)$ to denote that
$X$ has this probability law.  Let ${\rm supp}(X)$ denote the support of $X$.  Also, given an i.i.d.\ sequence $\{ X_i \}$, we will
often write $X$ for a ``generic" element of this sequence.

For any function $h$, let $h^\prime$ denote its first derivative,
and so on.  Also let ${\rm dom}\:( h)$ denote denote the domain of $h$
and ${\rm supp}\:( h)$ the support of $h$.  Finally, let ${\mathfrak
  F}_n$ denote the $\sigma$-field generated by the Letac's forward
sequence $V_0, V_1, \ldots, V_n$, where the process $\{ V_n \}$
is obtained from  the recursion \eqref{sr200}.
We now state the main hypotheses of this paper.\\[-.1cm]

\noindent
{\it Hypotheses}:\\
\indent
$(H_1)$  The random variable $A$ has an absolutely continuous component\\ 
\hspace*{1.25cm} with respect to Lebesgue measure, and satisfies the equation\\ 
\hspace*{1.25cm} $\Lambda(\xi) = 0$ for some  $\xi \in (0,\infty) \cap {\rm dom} \:( \Lambda^{\prime})$.

$(H_2)$   ${\bf E} \big[ | B |^\xi \big] < \infty$ and ${\bf E} \big[ \left(A|D| \right)^\xi \big] < \infty$.

$(H_3)$  ${\bf P} \left( A> 1, B> 0 \right) > 0$
or ${\bf P} \left( A>1, B \ge 0, D > 0 \right)>0$.\\[-.1cm]

Next we turn to a characterization for the constant of decay.  To describe this constant, we will compare the process $\{ V_n \}$
generated by Letac's forward equation \eqref{letac2} to that of a random walk, namely 
\[
S_n := \log A_1 + \cdots + \log A_n, \quad n=1,2,\ldots,
\]
where $\{ \log A_i \}$ is an i.i.d.\ sequence with common distribution $\mu_A$.
In Sparre-Andersen's generalization of the classical ruin problem,
one studies ruin for the random walk $\{ S_n \}$, that is,
\[
\Psi^\ast(u) := {\bf P} \left( S_n < -u, \:\: \mbox{\rm for some} \: \: n \in \pintegers \right),
\]
and it is well known that, asymptotically,
\begin{equation} \label{sr202}
\Psi^\ast(u) \sim C^\ast e^{-\xi u} \quad \mbox{\rm as} \quad u \to \infty,
\end{equation}
where ${\bf E} \big[ A^\xi \big] = 1$.
The constant $C^\ast$ may be described in various ways, but a characterization due to \shortciteN{DI72} is obtained
by first setting
\begin{equation} \label{sr202a}
\tau^\ast = \inf \left\{ n \ge 1:  S_n \le 0 \right\},
\end{equation}
and then defining
\begin{equation} \label{sr203}
C^\ast =  \frac{1- {\bf E} \big[ e^{\xi S_{\tau^{\star}}}\big]}{\xi \lambda^\prime(\xi) {\bf E} \left[ \tau^\ast \right]}.
\end{equation}
We will refer to $C^\ast$ as the Cram\'{e}r-Lundberg constant.

Our aim is to develop an analogous characterization for the constant of 
tail decay for the sequence $\{ V_n \}$.  Expectedly, it
will involve the return times $\tau$, where $\tau$ is a typical
regeneration time of the Markov chain in the sense of Lemma 2.2.
In Lemma \ref{MarkovReg} below, we will show that $\{ V_n \}$ is indeed a Markov chain satisfying appropriate path properties for this lemma to hold.

Returning to the central result of the paper, recall that
$\{Z_n^{(p)} \}$ and $\{ Z_n^{(c)} \}$ are the associated perpetuity and
conjugate sequences, respectively, and let  $Z^{(p)}$
and $Z^{(c)}$ denote their a.s.\ limits.  In Lemma 5.5 below, 
we will show that these limits 
exist and that $Z^{(p)} > Z^{(c)}$ a.s.  

%
%
%  THEOREM 2.1
%
%

\begin{thm} Assume Letac's Model E, and suppose that $(H_1)$, $(H_2)$, and $(H_3)$ are satisfied.  Then
\beq{\label{eq2.1}}
\lim_{u \ra \ff} u^{\xi}{\bf P} \left( V>u \right) =C
\eeq
for a finite positive constant $C$ which is given by
\begin{equation} \label{eq2.2}
C = \frac{1}{\xi \lambda^\prime(\xi) {\bf E} \left[ \tau \right]} {\bf
  E}_\xi \left[ \left( Z^{(p)} - Z^{(c)} 
\right)^\xi {\bf 1}_{\{ \tau = \infty \}} \right].
\end{equation}
\end{thm}

\begin{Rk}{\rm
Comparing \eqref{eq2.2} with \eqref{sr203},
we obtain alternatively that $C = C^\ast \theta  {\bf E}_\xi \big[ \big( Z^{(p)} - Z^{(c)} 
\big)^\xi {\bf 1}_{\{ \tau = \infty \}} \big]$, where
\[
\theta := \left( {\bf E} \left[ \tau \right] \right)^{-1} \frac{ {\bf E} \left[ \tau^\ast \right] }{ 1 - {\bf E} \left[ e^{\xi S_{\tau^\ast}} \right] }.
\]
}
\end{Rk}

It should be noted that the expectation in the previous theorem is
computed {\em in the} $\xi-${\em shifted measure}, and for this
reason, the event $\{ \tau = \infty \}$ occurs with positive probability.
Moreover, the constant $C$ obtained in Theorem 2.1 can be decomposed into three parts, of which we would like to emphasize the
{\it Cram\'{e}r-Lundberg constant} $C^\ast$ and ${\bf E}_\xi \big[ \big( Z^{(p)} - Z^{(c)} 
\big)^\xi {\bf 1}_{\{ \tau = \infty \}}  \big]$.
In the proof, it will be seen that the constant $C^\ast$ results from the ``large-time" behavior of the process which,
under the rare excursion to the high level $u$, resembles the random walk $S_n = \sum_{i=1}^n \log A_i$.
Indeed, when $(B,D)= (0, 1)$, then the constant $C$ reduces to the Cram\'{e}r-Lundberg constant $ C^\ast$.
In contrast, the constant ${\bf E}_\xi \big[ \big( Z^{(p)} - Z^{(c)} 
\big)^\xi {\bf 1}_{\{ \tau = \infty \}}  \big]$ results from the 
``small time" behavior of the process.   

For the elementary recursion defined by  $f(v) = Av + B$,  \shortciteN{NECSOZ09} also provide a probabilistic representation for $C$ but their characterization is complex and is not easily amenable for statistical computation.
In contrast,  we obtain a comparatively simple 
form for the constant, namely ${\bf E}_\xi \big[\big( Z^{(p)} - Z^{(c)} 
\big)^\xi {\bf 1}_{\{
  \tau = \infty \}}  \big]$,
although the exact form of this constant will depend on the particular
recursion under study.  However, in a wide variety of applications,
the conjugate term $Z^{(c)}$ will actually be zero on $\{ \tau = \infty \}$, and so the integrand
involves only the perpetuity term, $Z^{(p)}$, which has the simple
form of \eqref{aps3} and whose terms can be shown to converge to zero at a
geometric rate.

Returning to the ruin problem with investments described in Section 1.2, the forward recursive sequence 
is given by
\begin{equation} \label{sr210}
V_n = \left( A_n V_{n-1} + B_n \right)^+, \quad n = 1,2,\ldots, \quad
V_0 =v.
\end{equation}
First, to determine the regeneration time $\tau$ for this process, begin by observing that the minorization condition $({\cal M})$ holds
with
\begin{equation}\label{minor-ruin}
\delta = 1, \quad {\cal C} = \{ 0 \}, \quad \mbox{\rm and} \quad \nu = {\mathfrak L}( B^+).
\end{equation}
In other words, $\{ V_n \}$ has an atom at $\{ 0 \}$.
Consequently, by Remark 2.1, we may take $\tau -1$ to be the first return
time of the process $\{ V_n \}$ to the origin, starting from an
initial state $V_0$, where $V_0 \sim \nu =  {\mathfrak L}( B^+).$
Next we turn to identifying the conjugate sequence.
Upon setting $D_i = -B_i/A_i$
in \eqref{def-acs4}, we obtain
\begin{equation} \label{def-acs4a}
Z^{(c)} = \min \left( 0, \bigwedge_{n=1}^\infty \sum_{i=0}^n
  \frac{B_i}{A_0 \cdots A_i} \right).
\end{equation}
But the terms on the right-hand side all have the form
$\sum_{i=0}^n B_i/(A_0 \cdots A_i)$ for some $n$, which
upon multiplication by $(A_0 \cdots A_n)$ yields
\begin{equation} \label{ruin-suppl1}
\sum_{i=0}^n B_i \left( A_{i+1} \cdots A_n \right), \quad n=0,1,\ldots.
\end{equation}
Now for $Z^{(c)}$ to be nonzero, we would need the quantity
in \eqref{ruin-suppl1} to be negative for some $n$.  But
\eqref{ruin-suppl1} is just the $n$th iteration of the SFPE $\tilde{V}
\stackrel{\cal D}{=} A\tilde{V}+B$ with forward recursive sequence $\{
\tilde{V}_n \}$, and from the definition of $\{ V_n \}$ and $\{
\tilde{V}_n \}$, we have
\[
\inf \{n:  \tilde{V}_n \le 0 \} = \inf \{n:  V_n \le 0 \} := \tau -1,
\]
since these two processes {\em agree} up until their first entry into
$(-\infty,0]$,
at which time the reflected process $\{ V_n \}$ is then set equal to zero.
Hence we conclude that, conditional
on $\{ \tau = \infty \}$, $Z^{(c)} = 0$ and so the
constant $C$ is determined entirely by the perpetuity term $Z^{(p)}$,
and conditioning on $\{ \tau = \infty \}$ simply removes those paths
for which the process $\{ V_n \}$ returns to the origin in its
$\xi$-shifted measure.

Next we supplement the previous estimate with a Lundberg-type upper bound.  First 
define the perpetuity sequence
\begin{equation} \label{def-barW}
\bar{Z}^{(p)} = |V_0| +  \sum_{i=1}^\infty \frac{  \big( |B_i|  + A_i
  |D_i|  \big)  }{A_1 \cdots A_i}
  {\bf 1}_{\{ \tau > i \}}.
\end{equation}

%
%
%  THEOREM 2.2
%
%

\begin{thm} Assume the conditions of the previous theorem.  
Then
\beq{\label{sr95}}
{\bf P}\left(V>u \right) \le \bar{C}(u) u^{-\xi}, \quad \mbox{\rm for all}\:\:u \ge 0,
\eeq
where, for certain positive constants $C_1(u)$ and $C_2(u)$,
\begin{equation} \label{sr96}
\bar{C}(u) =  {\bf E}_\xi \big[\big( \bar{Z}^{(p)}\big)^\xi \big]  \sup_{z \ge 0} \Big\{ e^{-\xi z}\left( zC_1(u) + C_2(u) \right)\! \Big\}
< \infty \quad \mbox{\rm for all}\: u.
\end{equation}
Furthermore,
\begin{equation} \label{sr97}
C_1(u) \to \frac{1}{m} \quad \mbox{\rm and} \quad C_2(u) \to 1 + \frac{\sigma^2}{m^2} \quad \mbox{\rm as} \:\: u \to \infty,
\end{equation}
where $m := {\bf E}_\xi \left[ \log A \right]$ and $\sigma^2 := {\rm Var}_\xi \left( \log A \right)$.
\end{thm}

The constants $C_1(u)$ and $C_2(u)$ will be identified
explicitly in \eqref{pf-UB6} below.

\subsection{Generalizations of Letac's Model E} 
While the results of the previous section were limited to Letac's Model E, the method of proof is actually more general,
applying to a wide class of recursions.  For example, by a slight modification of the proof, it is also possible to 
analyze polynomial recursions of the form
\begin{equation} 
V \stackrel{\cal D}{=} B_k V + B_{k-1} V^{(k-1)/k} + \cdots + B_{1} V^{1/k} + B_0,
\end{equation}
where the random variable $B_k$ is assumed to satisfy the same conditions as the random variable $A$ in Letac's model,
and the vector $(B_0,\ldots,B_k)$ takes values in $(0,\infty)^{k+1}$.  We obtain the same asymptotic results as before, although the form of the
constant $C$ is more complicated.  We will return to this example and
 some further generalizations in Section 7, and
to an extension to a Markovian driving sequence in Section 8.

%
%
%  EXAMPLES AND APPLICATIONS
%
%
\section {Examples and applications}
In this section, we develop a few examples and draw a connection
to the related problem of estimating the extremal index.

\subsection{The ruin problem with stochastic investments}  Returning to the ruin problem described in Section 1.2, we recall that in 
this case the SFPE specializes to 
$f(v) = \left( Av+B\right)^+$.

In the context of this problem, it is interesting to point out that
there is a close connection of our approach to methods based on duality.  Specifically,
to obtain the above SFPE, we began by studying the risk
process in \eqref{intro2},
which leads,
via Letac's forward recursion, to the Markov process
\begin{equation} \label{ex1-new3}
V_n = \left( A_n V_{n-1} + B_n \right)^+, \quad n=1,2,\ldots, \quad V_0=0,
\end{equation}
which we will analyze in the proofs of our theorems below.
Alternatively, we could have arrived at the process in
\eqref{ex1-new3} using a duality technique as
described in \shortciteN{SAKS96}.
In particular, the main result of that article states that the
``content process'' $\{ V_n \}$ is dual to the ``risk process''
$\{Y_n \}$.  In particular, their risk process is obtained by inverting
the function $f(v;a,b) = (av + b)^+$ for fixed $a,b$, and then
studying the forward recursion generated by $f^{\leftarrow}$, which is
easily computed to be the function $\{ Y_n \}$ in \eqref{intro2}.
Hence by \shortciteN{SAKS96}, Theorem 3.1,
\begin{equation} \label{ex1-new4}
{\bf P} \left( V > u \right) = {\bf P} \left( Y_n < 0, \:\: \mbox{for some}\:\: n \in \pintegers \right) := \Psi(u).
\end{equation}
This shows---by an alternative approach---that the probability of ruin
may be equated to the steady-state exceedance probability of the
forward recursive sequence $\{ V_n \}$.

Appealing now to Theorem 2.1 of this paper, we obtain \eqref{eq2.1} with
\begin{equation} \label{ex1-11}
C = \frac{1}{Ê\xi \lambda^\prime(\xi) {\bf E} \left[ \tau \right]}
{\bf E}_\xi \left[ \left( \sum_{i=0}^{\infty}
 \frac{B_i}{A_0 \cdots A_i} \right)^\xi {\bf 1}_{\{ \tau = \infty\}}\right]
\end{equation}
where
$\tau - 1$ is the first return time of $\{ V_n \}$ to the origin 
and $B_0 \sim {\mathfrak L}(B^+)$; {\it
  cf.}\ the discussion following \eqref{minor-ruin}.
Note that the last expectation may be simplified for $\xi \ge 1$,
since Minkowskii's inequality followed by a  change
of measure argument yields the upper bound
\[
 {\bf E}_\xi \left[ \left( \sum_{i=0}^{\infty} \frac{|B_i|}{A_0 \cdots
       A_i} \right)^\xi {\bf 1}_{\{ \tau = \infty\}}\right]^{1/\xi}
 \le {\bf E} \left[ | B|^\xi \right]^{1/\xi} \left(1+ \sum_{i=1}^\infty  {\bf P} \Big( \tau > i-1 \Big)^{1/\xi}\right).
\]
Thus the constant $C$ may be explicitly related to the return times of the process $\{ V_n \}$ to the origin.  
Since $\{ V_n \}$ is geometrically ergodic (Lemma 5.1 below), these
regeneration times have exponential moments,
and so this last expression is also finite.

Finally, if $\xi < 1$, then an analogous
result holds by applying the deterministic inequality
\begin{equation} \label{ex1-13}
\left| x + y \right|^\xi \le \left|x\right|^\xi + \left|y \right|^\xi, \quad \forall (x,y) \in \reals^2,\: \xi \in (0,1], 
\end{equation}
in place of Minkowskii's inequality.
For some specific examples of investment and insurance processes, see \shortciteN{JC09}, Section 3. 

\subsection{The ARCH (1) and GARCH(1,1) financial processes}
From a mathematical perspective, the ruin problem described in the previous example is closely related to a problem in 
financial time series modeling, namely, the characterization of the
extremes of the GARCH(1,1) and ARCH(1) financial
processes.

In particular, in the GARCH(1,1) model of \shortciteN{TB86}, the logarithmic returns on an asset are modeled as
\begin{equation} \label{ex2-1}
R_n = \sigma_n \zeta_n, \quad n =1,2,\ldots,
\end{equation}
where $\{ \zeta_n \}$ is an i.i.d.\ Gaussian sequence of random variables, and where $\{ \sigma_n \}$ satisfies the recurrence equation
\begin{equation} \label{ex2-2}
\sigma_n^2 = a + b \sigma_{n-1}^2 + c R_{n-1}^2
\end{equation}
for certain positive constants $a$, $b$, and $c$.
The quantity $\sigma_n$ is called the stochastic volatility.  Substituting \eqref{ex2-1} into \eqref{ex2-2} yields 
the recurrence equation
\begin{equation} \label{ex2-3}
V_n = A_n V_{n-1} + B_n, \quad n=1,2,\ldots.
\end{equation}
with $V_n = \sigma_n^2$, $A_n = \left(b + c\xi_{n-1}^2 \right)$, and
$B_n = a$.  We note that this same recurrence equation also arises in the
well-known ARCH(1) model of \shortciteN{RE82}.

By comparing this SFPE with the SFPE of the previous example,
we see that we are dealing with essentially the same process as before, except that in this case the random variables $B_n$ are always positive,
meaning that we may no longer take $\tau$ to be the first return time to $\{ 0 \}$.   To determine a correct choice of $\tau$
we need to verify, as before, the minorization condition $({\cal M})$
of Section 2.1.   To this end, note that in the GARCH(1,1) model, $V(v):= Av+B$ has a density, which we call $h_v$.
Since $h_v(x)$ is monotonically decreasing for large $x$, 
\[
\inf_{v \in [0,M]} h_v(x) = h_0(x) \quad \mbox{\rm for sufficiently large}\:\: x.
\]
Hence 
\[
\hat{\nu}(E) := \int_E \left(\inf_{v \in [0,M]} h_v(x)\right) dx
\]
is not identically zero, and then 
$P(x,E)Ê\ge {\bf 1}_{[0,M]}(x) \hat{\nu}(E).$  Finally, upon normalizing $\hat{\nu}$
so that it is a probability measure, we obtain $({\cal M})$ with $k=1$ and $\delta$ the normalizing constant.
Thus, according to  Remark 2.1, we see that regeneration occurs w.p.\ $\delta$ upon return to the set ${\cal C} = [0,M]$.

Once this modification has been made in the choice of $\tau$, 
we may repeat the calculation in \eqref{ex1-11} of the previous example to obtain an
explicit characterization of the constant $C$ in Theorem 2.1.  As in
the previous example, we obtain that the conjugate term $Z^{(c)}$ is
zero on $\{ \tau = \infty \}$, in this case because the process $\{
V_n \}$ is {\em always} positive, and so we indeed obtain
\eqref{ex1-11} but with a different $\tau$.

\subsection{The extremal index}  
While the previous examples and estimates have dealt with stationary tail probabilities, the methods can also be applied
to other related problems.  The extremal index---which measures the tendency of
a dependent process to cluster---is defined for a strictly stationary process $\{ \tilde{V}_n \}$ by the equation
\begin{equation} \label{index1}
{\bf P} \left( \max_{1 \le i \le n} \tilde{V}_n \le u_n \right) = e^{-\Theta t},
\end{equation}
for any given $t$, where $\{ u_n \}$ is chosen such that
\begin{equation} \label{index2}
\lim_{n \to \infty} n \left\{ 1 - {\bf P} \left( \tilde{V}_n > u_n \right) \right\} = t.
\end{equation}
The constant $\Theta$ in \eqref{index1} is called the extremal index.
As argued in \shortciteN{HR88}, p.\ 380,
this quantity may be written in a Markovian setting as
\begin{equation} \label{ex3-1}
\Theta = \lim_{u \to \infty} \frac{ {\bf P} \left( V_n > u, \: \mbox{\rm some}\: n < \tau \right)}{{\bf E} \left[ N_u \right] },
\end{equation} 
where $N_u$ denotes the number of exceedances above level $u$ occurring over a regeneration
cycle.  We will show in this article, using the methods
of Theorem 2.1, that
\begin{equation} \label{ex3-6}
\Theta = \frac{1 - {\bf E} \big[ e^{\xi S_{\tau^\ast}} \big]}{{\bf E} \left[ \tau^\ast \right] },
\end{equation}
where $\tau^\ast = \inf \{ n \ge 1:  S_n \le 0\}$ for $S_n = \sum_{i=1}^n \log A_i$.
The above expression provides a considerable simplification of known theoretical estimates for $\Theta$, as given, {\it e.g.}, 
for the special case of the GARCH(1,1) model in
\shortciteN{LHHRSRCV89} and, moreover, this last expression is seen to hold in the general setting of Letac's Model E.

%
%
%  SOME RESULTS FROM NONLINEAR RENEWAL THEORY
%
%
\section{Some results from nonlinear renewal theory}
The main tools needed in the proofs of Theorems 2.1 and
2.2 will involve ideas from nonlinear renewal theory,
which we now present in detail.

Recall that the crux of the proof of Theorem 2.1 will be to utilize the representation formula of
Lemma 2.2, namely,
\[
{\bf P} \left( V > u \right) = \frac{ {\bf E} \left[  N_u \right]}{{\bf E} \left[ \tau \right]}.
\]
Thus we will need to estimate ${\bf E} \left[  N_u \right]$, the number of exceedances above level $u$ which occur {\em over a regeneration cycle}.
To study this quantity, we will employ a  {\em dual} change of measure, defined as follows.  First let $(A_1,B_1,D_1) \sim \mu_\xi$,
where $\mu_\xi$ is the $\xi$-shifted measure defined previously in \eqref{meas-shift}.  If $V_1 > u$, then let
$(A_n,B_n,D_n) \sim \mu$ for all $n > 1$.  Otherwise let $(A_2,B_2,D_2) \sim \mu_\xi$; and then for $V_2 > u$ set $(A_n,B_n,D_n) \sim \mu$ for all $n > 2$,
while for $V_2 \le u$ set $(A_3,B_3,D_3) \sim \mu_\xi$, etc.
In summary, the dual change of measure can be described as follows:
\begin{equation} \label{DUAL-MEAS}
{\mathfrak L} \big(\log A_n,B_n,D_n\big) = \left\{ \begin{array}{l@{\quad}l}
\mu_\xi & \mbox{for}\:\:  n=1,\ldots,T_u, \\
\mu & \mbox{for}\:\: n > T_u,
\end{array} \right.
\end{equation}
where $T_u = \inf \{n:  V_n > u \}$ and $\mu_\xi$ is defined as in \eqref{meas-shift}.
Roughly speaking, this dual measure shifts the distribution of $\log A_n$
on a path terminating at time $T_u$, and reverts to the original measure thereafter. 
Let ${\bf E}_{\mathfrak D} \left[ \cdot Ê\right]$ denote expectation with respect to the dual measure described in \eqref{DUAL-MEAS}. 

We now focus on two quantities:  
\begin{description}
\item (i) the overjump distribution at the first time $V_n$ exceeds
the level $u$,\\
\hspace*{-.3cm} calculated in the $\xi$-shifted measure; and 
\item (ii) the expected number of exceedances above level $u$ which then occur
up\\
\hspace*{-.25cm} until regeneration, calculated in the original measure.  
\end{description}
Note that in the $\xi$-shifted measure,
${\bf E}_\xi \left[ \log A \right] > 0$ and hence $V_n \uparrow
\infty$ w.p.1 as $n \to \infty$ ({\it cf.}\ Lemma 5.2 below).  Consequently,
\[
V_n = A_n \max \left( D_n, V_{n-1} \right) + B_n
\]
implies 
\[
V_n \approx A_n V_{n-1} \quad \mbox{\rm for large}\:\: n.
\]
In other words, for large $u$ the process $\{ V_n \}$ will ultimately resemble a perturbation of multiplicative random walk in an asymptotic sense.
Consequently, the problem described in (i) may be viewed as a variant of results from nonlinear renewal theory
({\it cf.}\ \shortciteN{MW82}).  

%
%
%  LEMMA 4.1
%
%

\begin{Lm}
Assume Letac's Model E, and suppose that $(H_1)$, $(H_2)$, and $(H_3)$ are satisfied.  Then in the dual measure \eqref{DUAL-MEAS},
\begin{equation} \label{prelimNR1}
\lim_{u \to \infty} {\bf P}_\xi \left( \left. \frac{V_{T_u}}{u} > y \right| T_u <
  \tau \right) = {\bf P}_\xi \left( \hat{V} > y \right)
\end{equation}
for some random variable $\hat{V}$.  The distribution of this random variable $\hat{V}$ is independent of the initial distribution of $V_0$
and is described as follows.
If $A^l$ is a typical ladder height of the process $S_n = \sum_{i=1}^n \log A_i$
in the $\xi$-shifted measure, then
\begin{equation} \label{prelimNR2}
{\bf P}_\xi \big( \log \hat{V} > y \big) = \frac{1}{{\bf E}_\xi \big[ A^l \big]} \int_y^\infty {\bf P}_\xi \big(A^l > z \big) dz, \quad \mbox{\rm for all}\:\: y \ge 0.
\end{equation}
\end{Lm}

While Lemma 4.1 follows rather immediately from known results in nonlinear renewal theory, this is not the case
for the problem stated in (ii) above.  Here we would like to determine
${\bf E}_{\mathfrak D} \left[ N_u \left| {\mathfrak F}_{T_u} \right. \right]$, {\it i.e.}, the expected number of exceedances, in the dual measure, which occur
over a typical regeneration cycle prior to
the regeneration time $\tau$.
We provide a precise estimate for this quantity next.
In the following result, note that we {\em start} the process at a level
$vu$, where $v>1$, and so the dual measure actually agrees with
the original measure in this case. 

%
%
%  THEOREM 4.1
%
%

\begin{thm}
Assume Letac's Model E, and suppose that $(H_1)$, $(H_2)$, and $(H_3)$ are satisfied.  Then for any $v > 1$,
\beq
\lim_{u \ra \ff} {\bf E} \left[  N_u \, \left| \, \frac{V_0}{u} = v \right. \right] = U(\log v),
\eeq
where  $U(z):= \sum_{n \in {\mathbb N}} \mu_A^{\star n}(-\ff, z)$  and $\mu_A$ is the marginal distribution of $-\log A$.  
\end{thm}

Roughly speaking, the function $U$ may be interpreted as the renewal function of the random walk
$-S_n = - \sum_{i=1}^n \log A_i$.  More precisely, if the distribution of $\log A$ is continuous, so that the open interval $(-\infty,z)$ could
be replaced with the closed interval $(-\infty,z]$, then $U(z)$ agrees with the standard definition of the renewal function of $\{-S_n \}$.

We now supplement the estimate in the previous theorem with an upper bound.

%
%
%  THEOREM 4.2
%
%

\begin{thm}Assume that the conditions of the previous theorem.  Then there exist finite positive constants $C_1 (u)$ and $C_2(u)$ such that
\begin{equation} \label{nr10}
{\bf E}_{\mathfrak D} \left[ N_u \, \left| \, {\mathfrak F}_{{T_u}\wedge (\tau-1)}\right. \right] \le \left(C_1(u) \log \left( \frac{V_{T_u}}{u} \right) 
+ C_2(u)\right) {\bf 1}_{\{T_u < \tau \}},
\:\: \mbox{\rm all } u,
\end{equation}
where $C_1(u)$ and $C_2(u)$ are the positive finite constants
that converge as $u \to \infty$ to $m^{-1}$ and $1+(\sigma^2/m^2)$, respectively,
where $m := {\bf E}_\xi \left[ \log A \right]$ and $\sigma^2 := {\rm Var}_\xi \left( \log A \right)$.
\end{thm}
We emphasize that this last theorem is crucial for obtaining the Lundberg upper bound in Theorem 2.2.

Finally, using the previous results, it is now possible to state an extension of Lemma 4.1, which will be particularly useful
in the sequel.  First let 
\[
{\cal Q}_u := {\bf E}_{\mathfrak D} \left[ N_u
  \left| {\mathfrak F}_{T_u \wedge (\tau-1)} \right. \right] \left(
  \frac{V_{T_u}}{u} \right)^{-\xi}
{\bf 1}_{\{ T_u < \tau \}}.
\]

In the following, we write ${\bf
  P}_{\mathfrak D}(A)$ for ${\bf E}_{\mathfrak D} \left[ {\bf 1}_A \right]$.

%
%
%  THEOREM 4.3
%
%

\begin{thm}
Assume the conditions of the previous theorem.  Then conditional on $\{ T_u < \tau \}$,  
\begin{equation} \label{pf-oldCor}
{\cal Q}_u \Rightarrow U(\log \hat{V}) \hat{V}^{-\xi} \quad \mbox{\rm as} \quad u \to \infty,
\end{equation}
where $\hat{V}$ is given as in \eqref{prelimNR1} and \eqref{prelimNR2}.  That
is, 
\[
\lim_{u \to \infty} {\bf P}_{\mathfrak D} \left(\left. {\cal Q}_u \le y \right| T_u < \tau
  \right) = {\bf P}_{\mathfrak D} \left( U(\log \hat{V})
  \hat{V}^{-\xi} \le y\right),
\quad \mbox{\rm for all} \:\: y \ge 0.
\]
\end{thm}

%
%
%  PROOFS OF MAIN RESULTS
%
%
\section{Proofs of the main results}
\subsection{Five preliminary lemmas}
We begin by presenting five preparatory lemmas, which will be needed in the proofs of the main theorems to be
proved in Section 5.2.   The first two lemmas establish
path properties of the Markov chain $\{ V_n \}$, obtained via 
forward iteration of the SFPE.
Recall that a Markov chain satisfies a {\em drift condition} if
\begin{equation}
\int_{\mathbb S} h(y) P(x,dy) \le \rho h(x) + \beta {\bf 1}_{\cal C}(x), \quad \mbox{\rm for some}\:\: \rho \in (0,1),
 \tag{${\cal D}$}
\end{equation}
where $h$ is a function taking values in $[1,\infty)$, $\beta$ is a
positive constant, and  ${\cal C}$ is a 
Borel subset of $\reals$.  In the next theorem, we work with the assumption that the process
$\{ V_n \}$ is nondegenerate, namely that we do not have ${\bf P} \left( F_Y(v) = v \right)=1$.
This assumption is subsumed by $(H_3)$, but we state it here as a
separate hypothesis.

%
%
%  LEMMA 5.1
%
%

\begin{Lm} \label{MarkovReg}  Assume Letac's Model E, and let $\{ V_n \}$ denote the forward recursive sequence
corresponding to this SFPE.
Assume that $(H_1)$ and $(H_2)$ are satisfied and that the process $\{
V_n \}$ is nondegenerate.   Then{\rm :}

{\rm (i)}  $\{ V_n \}$ satisfies the drift condition $({\cal D})$ with ${\cal C} = [-M,M]$.

{\rm (ii)}  $\{ V_n \}$ is $\varphi$-irreducible, where $\varphi$ is
the stationary distribution of $\{ V_n \}$.

{\rm (iii)}   $\{ V_n \}$ satisfies the minorization condition $({\cal M})$ with $k=1$ under both the measure $\mu$ and  the measure $\mu_\xi$.  Furthermore, the set $[-M,M]$ is petite for any $M > 0$.

{\rm (iv)}  $\{ V_n \}$ is geometrically ergodic.  Moreover
${\bf E} \left[ e^{\epsilon \tau} \right] < \infty,$
where $\tau$ is the inter-regeneration time of the process $\{ V_n \}$
under {\rm any} minorization $({\cal M})$
where the ${\cal C}$-set is bounded.
\end{Lm}  

\Pf
(i)   Using the defining equation 
$V_n = A_n \max \left( D_n, V_{n-1} \right) + B_n$, we obtain the inequality
\begin{equation} \label{pf-drift1}
|V_n| \le A_n |V_{n-1}| + \left( A_n |D_n| + |B_n| \right).
\end{equation}
Now hypothesis $(H_1)$ implies that ${\bf E} \left[ A^\alpha \right] < 1$ for all $\alpha \in (0,\xi)$.  Fix $\alpha \in \left( \xi \wedge 1 \right)$
and $V_0 = v$.
Then from the deterministic inequality \eqref{ex1-13}, we obtain
\begin{equation} \label{Markovreg-1}
{\bf E} \left[ | V_1|^\alpha \right] \le {\bf E} \left[ A^\alpha \right] v^\alpha + {\bf E} \left[ |B_1|^\alpha \right] + {\bf E} \big[ \left(A_1 |D_1| \right)^\alpha \big].
\end{equation}
Thus $({\cal D})$ holds with $h(x) = |x|^\alpha +1$, $\rho = \left( {\bf E} \left[ A^\alpha \right] + 1 \right)/2$, and ${\cal C} = [-M,M]$, where 
$M$ is chosen sufficiently large such that
\[
\frac{1}{2} \left( 1 - {\bf E} \left[ A^\alpha \right]  \right)M^\alpha \ge {\bf E} \big[ |B_1|^\alpha \big] + {\bf E} \big[ (A_1 |D_1|)^\alpha \big].
\]

(ii)  To verify that $\{ V_n \}$ is irreducible, it is sufficient to
show that if $\varphi(A) > 0$, then ${\cal U}(x,A) > 0$ for all $x$,
where
\[ 
{\cal U}(x,A) := \sum_{n=1}^\infty P^n(x,A)
\]
({\it cf.}\ \shortciteN{SMRT93}, Proposition 4.2.1).  Now to identify
$\varphi$, note under our hypotheses that the conditions of
Proposition 6.1 of \shortciteN{CG91} are satisfied. Hence the backward
recursion \eqref{eq0.15} converges and is independent of $v$.  Thus,
by Lemma 2.1, it follows that $\{ V_n \}$ is stationary and the law of
the limiting random variable $V$ is the same as the law of $Z$.  Now
choose $\varphi \equiv {\mathfrak L}(Z)$.  Since $\{ V_n \}$ is
aperiodic, it follows by Theorem 1.3.1 and (1.16) of
\shortciteN{SMRT93}) that
\begin{equation}
\sup_{A \in {\cal B}(\reals)} \left| P^n(x,A) - \varphi(A)
\right| \to 0 \quad \mbox{\rm as}\:\: n \to \infty,\quad \forall x \in
\reals
\end{equation}
(where we have used (i) to verify their condition (1.14) with
$h(x) = |x|^\alpha + 1$).  Hence
\begin{equation}
{\cal U}(x,A) \ge  \lim_{n \to \infty} P^n(x,A) = \varphi(A)
> 0
\end{equation}
for any $\varphi$-positive set $A$.

(iii)  We note that the main issue is to verify that the minorization
holds with $k=1$ (rather than for general $k$).

We will show that for any $v \in \reals$,
there exists an $\epsilon$-neighborhood ${\mathfrak B}_\epsilon(v)$ such that
\begin{equation} \label{minor-1}
P(v,E) \ge \delta {\bf 1}_{{\mathfrak B}_\epsilon(v)} \nu(E), \quad \mbox{\rm for all}\:\: E \in {\cal B}(\reals),
\end{equation}
for some positive constant $\delta$ and some probability measure 
$\nu$, both of which will typically depend on $v$.  Since {\em some}\,
such interval ${\mathfrak B}_\epsilon(v)$ will necessarily be $\varphi$-positive
with $\varphi = {\mathfrak L}(V)$ as in (ii), this will imply the existence of a minorization with $k=1$.

Set  $v^\ast = \inf
\left\{ v:  {\bf P} \left( D \le v \right)=1 \right\} \in (-\infty,\infty]$.
We will consider three different cases, namely $v < v^\ast$, $v > v^\ast$,
and $v=v^\ast$.

If $v < v^\ast$, then ${\bf P} \left( D > v \right) > 0$.  Since ${\bf P} \left( D > w \right)$ is nonincreasing as a function 
of $w$, 
${\bf P} \left( D > v +\epsilon \right) > 0$ for some $\epsilon > 0$.
For this choice of $v$ and $\epsilon$, let ${\mathfrak B}_\epsilon(v)$
be the required $\epsilon$-neighborhood in \eqref{minor-1}.
Now if $V_0 \in {\mathfrak B}_\epsilon(v)$ and $D_1 > v + \epsilon$, then $\max( D_1,V_0 ) = D_1$ and hence
 $V_1 = A_1 D_1 + B_1$.   Thus, 
for any initial state $V_0 \in {\mathfrak B}_\epsilon(v)$,
\begin{equation} \label{minor-1a}
V_1 = \left( A_1 D_1 + B_1 \right) {\bf 1}_{\{ D_1 > v+\epsilon \}} + V_1 {\bf 1}_{\{ D_1 \le v + \epsilon \}}.
\end{equation}
Note that the first term on the right-hand side is independent of $V_0$.
Taking $\nu$ to be the probability law of $\left( A_1 D_1 + B_1 \right)$ conditional on $\{ D_1 > v+\epsilon \}$
and $\delta = {\bf P} \left( D > v +\epsilon \right)$,
we obtain \eqref{minor-1}.

Next suppose $v > v^\ast$.  Then for some $\epsilon > 0$, ${\bf P} \left( D > v - \epsilon \right) = 0$, and for
this choice of $v$ and $\epsilon$, let ${\mathfrak B}_\epsilon(v)$
be the $\epsilon$-neighborhood in \eqref{minor-1}.
Then ${\bf P} \left( V_1 = A_1 V_0 + B_1 \left|
 V_0 \in {\mathfrak B}_\epsilon(v) \right. \right) =1$.
Hence, to obtain a minorization for $V_1$ in this case,
it is sufficient to derive a minorization for $V_1 (v) := A_1 v + B_1$ a.s.
To this end, begin by observing that since $A_1$ has
 an absolutely continuous component, so does the pair $(A_1, A_1 v + B_1)$.
Let $\gamma^{(v)}$ denote the probability law of $(A_1, A_1 v + B_1)$.
Then by the Lebesgue decomposition theorem, $\gamma^{(v)}$
can be decomposed
into its absolutely continuous and singular components relative
 to Lebesgue measure, and the continuous component satisfies
\[
\gamma_c^{(v)} (E) = \int_E \frac{d\gamma^{(v)}}{dl}(z) dl(z), \quad \mbox{\rm for all}Ê\:\: E \in {\cal B}(\reals^2),
\]
where $l$ denotes Lebesgue measure on $\reals^2$. 
Hence by construction, there exists a rectangle where
$\frac{d\gamma^{(v)}}{dl}$ is bounded away from zero; that is,
\begin{equation} \label{minor-2}
\frac{d\gamma^{(v)}}{dl}(x,y) \ge \delta > 0, \quad \mbox{\rm for all}\:\: (x,y) \in [a,b] \times [c,d].
\end{equation}
Now suppose $w \in {\mathfrak B}_\epsilon(v)$ and let $\gamma^{(w)}$ denote the probability law of 
$V_1(w):= A_1 w + B_1$.  Then $V_1(w) - V_1(v) = A_1(w-v)$, and hence \eqref{minor-2} yields
\[
\frac{d\gamma^{(w)}}{dl}(x,y) \ge \delta, \quad \mbox{\rm for all} \quad (x,y) \in [a,b] \times [c+w^\ast,d-w^\ast],
\]
where $w^\ast = b|w-v|$.  
Notice that the constant $b$ was obtained by observing
that the first component describes the distribution of $A$, where
${\rm supp}(A) \subset (0,\infty)$, implying that $0 < a < b < \infty$.  Since $w^\ast \downarrow 0$ as $w \to v$, it follows that for sufficiently small
$\epsilon$, there exists a subinterval $[c^\prime,d^\prime]$ of $[c,d]$ such that
\begin{equation} \label{minor-3}
\frac{d\gamma^{(w)}}{dl}(x,y) \ge \delta, \quad \mbox{\rm for all} \:\: (x,y) \in [a,b] \times [c^\prime,d^\prime] \:\:
\mbox{\rm and all} \:\: w \in {\mathfrak B}_\epsilon(v).
\end{equation}
Hence the minorization \eqref{minor-1} holds with
\[
\nu(E) = \int_{\reals \times E} \: \left( \inf_{w \in {\mathfrak B}_\epsilon(v)} \frac{d\gamma^{(w)}}{dl}(x,y) \right) dl(x,y), \quad \mbox{\rm for all}\:\: E \in {\cal B}(\reals),
\]
and by \eqref{minor-2}, this measure is not identically equal to zero.

It remains to consider the case where $v = v^\ast$.  First observe that, similar to \eqref{minor-1a},
${\bf P} \left( V_1 = A_1 V_0 + B_1 \left| V_0 \in {\mathfrak B}_\epsilon
(v^\ast),\:\: D_1 \le v^\ast - \epsilon \right. \right) =1$.  Hence, we clearly have 
\begin{equation} \label{minor-9}
V_1 = \left( A_1 V_0 + B_1 \right) {\bf 1}_{\{ D_1 \le  v^\ast - \epsilon \}} + V_1 {\bf 1}_{\{ D_1 > v^\ast- \epsilon \}}.
\end{equation}
Now define $\gamma^{(\epsilon,v^\ast)}$ on any Borel set $E \subset \reals^2$ to be
\[
\gamma^{(\epsilon,v^\ast)}(E) = {\bf P} \left( (A_1, A_1v^\ast +B_1) 
  \in E , D_1 \le v^\ast-\epsilon
  \right) \nearrow \gamma^{(v^\ast)}(E) \quad \mbox{\rm as} \quad \epsilon
  \to 0,
\]
where the last step follows from the definition of $\gamma^{(v^\ast)}$, since
${\bf P} \left( D_1 \le v^\ast - \epsilon \right) \uparrow 1$ as
$\epsilon \to 0$.
Hence, since $\gamma^{(v^\ast)}$ has an absolutely continuous component,
there exists an $\epsilon > 0$ such that
$\gamma^{(\epsilon,v^\ast)}$ also has an absolutely continuous component.
Now apply the previous argument with $\gamma^{(v)}$ replaced with
$\gamma^{(\epsilon,v^\ast)}$ for this choice of $\epsilon$ to obtain
the corresponding minorization for this case.  Thus we have obtained
\eqref{minor-1} for all the three cases.
Note that the above computations hold regardless of whether we are in the original measure $\mu$ or in the $\xi$-shifted measure $\mu_\xi$.

Finally, to show that $[-M,M]$ is a petite set, 
note that $[-M,M] \subset \bigcup_{v \in [-M,M]} {\mathfrak B}_{\epsilon(v)}(v)$,
where ${\mathfrak B}_{\epsilon(v)}(v)$ is a small set and hence is petite.
Thus there exists a finite subcover of petite sets, and then by
Proposition 5.5.5 of \shortciteN{SMRT93}, $[-M,M]$ is petite.

(iv)  Since $[-M,M]$ is petite for any $M >0$, it follows from (i) and \shortciteN{SMRT93}, Theorem 15.0.1,
that the process $\{ V_n \}$ is geometrically ergodic.
It remains to show that, regardless of
$(\delta,{\cal C},\nu)$ in the minorization $({\cal M})$,
${\bf E} \left[ e^{\epsilon \tau} \right] < \infty$, where $\tau$ is
the inter-regeneration time under $({\cal M})$.  Now
since ${\cal C}$ is bounded, it follows by
Theorem 15.2.6 of
\shortciteN{SMRT93} that 
${\cal C}$ is $h$-geometrically regular
with $h(x) = |x|^\alpha +1$. Consequently, letting $K$
 denote the first return time of $\{ V_n \}$ to ${\cal C}$, we have
\begin{equation} \label{NEW-geoerg1}
\Gamma(t) := \sup_{v \in C} {\bf E} \left[ t^{K} | V_0 = v \right] < \infty
\end{equation}
for some $t > 1$.

Consider the split-chain (see
\shortciteN{EN84}, Section 4.4 and Remark 2.1 above).  Originating
from the ${\cal C}$-set, this chain has the 1-step transition measure
$\nu(dy)$ w.p.\:$\delta$ and the 1-step transition kernel
$\bar{P}(x,y):= (P(x,dy)-\nu(dy))/(1-\delta)$ w.p.\:$(1-\delta)$.  
Thus, conditional
on either one of these transition kernels, it follows from
\eqref{NEW-geoerg1} that $\Gamma_\nu(t)< \infty$ and 
$\Gamma_{\bar{P}}(t) < \infty$, where $\Gamma_\nu$ 
denotes conditioning on $V_1 \sim \nu$, and
$\Gamma_{\bar{P}}(t) :=  \sup_{v \in C} {\bf E} \left[ t^{K} | V_0 = v,
V_1 \sim \bar{P} \right]$.
Hence by Remark 2.1,
\begin{equation} \label{NEW-geoerg10}
{\bf E} \left[ t^\tau \right] \le \Gamma_\nu(t)
\left( \delta + \sum_{n=1}^\infty \delta (1-\delta)^{n} \left(\Gamma_{\bar{P}}(t)
  \right)^n  \right),
\end{equation}
where $\delta$ is the constant appearing in $({\cal M})$.
Finally observe by a dominated convergence argument that
$\Gamma_{\bar{P}}(t) \downarrow 1$ as $t \downarrow 1$, and thus for sufficiently small $t > 1$, 
$\Gamma_{\bar{P}}(t) < \left(1- \delta \right)^{-1}$. 
For this choice of $t$, the series on the right-hand side of 
\eqref{NEW-geoerg10} is convergent, as required.
\halmos

In the next result, we establish a critical result concerning the
transience of the process $\{ V_n \}$ in its $\xi$-shifted measure.

%
%
%  LEMMA 5.2
%
%

\begin{Lm}
Assume Letac's Model E, and let $\{ V_n \}$ denote the forward recursive sequence
corresponding to this SFPE.
Assume that $(H_1)$, $(H_2)$ and $(H_3)$ are satisfied.
Then under the measure $\mu_\xi$,
\begin{equation} \label{trans-1}
V_n \nearrow +\infty  \quad \mbox{\rm w.p.}\: 1 \quad\mbox{\rm as}\quad n \to \infty.
\end{equation}
Thus, in particular, the Markov chain $\{ V_n \}$ is transient.
\end{Lm} 

\Pf
We begin by showing that the chain $\{ V_n \}$ is transient
in the $\xi$-shifted measure.  
To this end, first note by $(H_3)$ that the set $[M,\infty)$ is attainable
with positive probability, for any positive constant $M$ and any initial state.  Thus $\varphi\left( [M,\infty) \right) > 0$ for all $M$.
Consequently by \shortciteN{SMRT93}, Theorem 8.3.6, it is sufficient to verify that for some positive constant $M$,
\begin{equation} \label{trans-2}
{\bf P}_\xi \left( V_n \le M, \: \mbox{\rm for some}\: n \in \pintegers \left| V_0 \ge 2M \right. \right) < 1.
\end{equation}

To establish \eqref{trans-2}, first note by definition of $\{ V_n \}$ that
\begin{equation} \label{trans-3}
V_n \ge A_n V_{n-1} - |B_n|, \quad \mbox{\rm for all} \: n,
\end{equation}
and iterating this equation yields
\begin{equation} \label{trans-4}
\frac{V_n}{A_1 \cdots A_n} \ge V_0 - W_n, \quad \mbox{\rm where} \quad W_n := \frac{|B_1|}{A_1} + \cdots + \frac{|B_n|}{A_1 \cdots A_n}.
\end{equation}

We will now show that $\{ W_n \}$ converges in distribution to a proper
nondegenerate random variable.  To this end, note that 
since $\{ (A_n,B_n ) \}_{n \in \pintegers}$ is an i.i.d.\ sequence, $W_n$ has the same distribution as 
\[
\tilde{V}_n := \sum_{i=1}^n \frac{|B_i|}{A_i} \prod_{j=i+1}^n A_j^{-1},
\]
where $\{ \tilde{V}_n \}$ is a Markov chain.  Moreover,
$\{ \tilde{V}_n \}$ satisfies hypotheses ($H_1$) and ($H_2$) and
is nondegenerate.  Consequently by Lemma 5.1 (iv), $\{ \tilde{V}_n \}$
converges to a proper distribution and hence so does $\{ W_n \}$.

Next observe that ${\bf E}_\xi \left[ \log A \right] = \Lambda^\prime(\xi) > 0$; thus,
 the i.i.d.\ sequence $\{ \log A_n \}$ has a positive drift in the $\xi-$shifted measure.  Consequently, setting $S_n = \sum_{i=1}^n \log A_i$
 we obtain ${\bf P} \left( S_n < 0, \:\: \mbox{\rm for some}\:\:n \right) < 1.$  Therefore
 \begin{equation} \label{trans-5}
 {\bf P}_\xi \left( A_1 \cdots A_n  \ge 1, \:\: \mbox{\rm for all}\:\:n \right) = p > 0.
 \end{equation}

Now assume that \eqref{trans-2} fails.  Then by \eqref{trans-4} and \eqref{trans-5}, we would have
\begin{equation} \label{trans-6}
{\bf P}_\xi \left( W_n \ge M, \:\: \mbox{\rm for some}\:\: n \right) \ge p, \quad \mbox{\rm for all} \:\: M.
\end{equation}
But $\{ W_n \}$ is nondecreasing and converges to a proper distribution.
Hence \eqref{trans-6} is impossible, so we conclude \eqref{trans-2}.

Next observe by Lemma 5.1 that the set $[-M,M]$ is petite for every $M >0$.
Hence $[-M,M]$ is uniformly transient
(\shortciteN{SMRT93}, Theorem 8.0.1) and, consequently,
the expected number of visits of $\{ V_n \}$ to $[-M,M]$ is finite for all $M$.
It follows that $| V_n | \uparrow + \infty$ w.p.\ 1 as $n \to \infty$.
But we cannot have $V_n \downarrow -\infty$
with positive probability since, by definition, we have as a lower bound that
$V_n \ge A_n D_n + B_n$ for each $n$, where $\{(B_n,D_n)\}$ forms an i.i.d.\ sequence of random variables. Consequently $V_n \uparrow + \infty$
w.p.1 as $n \to \infty$.
\halmos

Prior to stating the next lemma, recall that $T_u := \inf \left\{ n:  V_n > u \right\}$, and that ${\bf E}_{\mathfrak D}\left[ \cdot \right]$ denotes expectation with respect to the dual measure described in \eqref{DUAL-MEAS}.  We now compare an expectation in the original measure to this expecation in the dual measure,
following a typical regeneration cycle from its initial state with $V_0 \sim \nu$, where $\nu$ is the measure described in Lemma 2.2 (iv), and terminating
at the next regeneration time $\tau$.

%
%
%  LEMMA 5.3
%
%

\begin{Lm} 
Assume the conditions of Lemma 5.2.  Let $g:  \reals^\infty \to \reals$ be a deterministic function, and let $g_n$ denote its
projection onto the first $n+1$ coordinates; that is,
$g_n(x_0,x_1,\ldots) = g(x_0, \ldots, x_n)$.
Then
\begin{align} \label{measchange-1}
{\bf E} \left[ g_{\tau-1}(V_0,\ldots, V_{\tau-1}) \right] &= {\bf E}_{\mathfrak D}  \left[ g_{\tau-1}(V_0, \ldots, V_{\tau-1}) e^{-\xi S_{T_u}} {\bf 1}_{\{T_u < \tau\}}  \right]\\ \nonumber
& \;\; + {\bf E}_{\mathfrak D} \left[  g_{\tau-1}(V_0, \ldots, V_{\tau-1}) e^{-\xi S_{\tau}} {\bf 1}_{\{ T_u \ge \tau\}} \right].
\end{align}
\end{Lm}

\Pf
If
\[
{\mathfrak L} \left(\log A_k,B_k,D_k \right) = \left\{ \begin{array}{l@{\quad}l}
\mu_\xi & \mbox{for}\:\:  k=1,\ldots,n, \\
\mu & \mbox{for}\:\: k>n,
\end{array} \right.
\]
then it can be shown using induction that
\begin{equation} \label{pf-lm5.3.101}
{\bf E} \left[ g_n(V_0,\ldots,V_n) \right]
  = {\bf E}_\xi \left[g_n(V_0,\ldots,V_n)e^{-\xi S_n} \right].
\end{equation}
Now \eqref{measchange-1} follows by conditioning on $\{ T_u=m,\tau=n \}$
and summing over all possible values of $m$ and $n$.
\halmos

Next we establish a critical result which links Letac's forward
iteration of the SFPE to the backward equations.

%
%
%  LEMMA 5.4
%
%

\begin{Lm}
Assume Letac's Model E, and let $\{ V_n \}$ denote the forward recursive sequence
corresponding to this SFPE.  Let $\{ Z_n^{(p)} \}$ and $\{ Z_n^{(c)} \}$
denote the associated perpetuity sequence and the conjugate sequence,
respectively, assumed to have the initial values described above
in \eqref{def-aps1} and \eqref{def-acs1}, and set $A_0 = 1$.  Then for any
$n \in {\mathbb N}$,
\begin{equation} \label{prelm3-1} 
\left(Z_n^{(p)} - Z_n^{(c)} \right){\bf 1}_{\{ \tilde{Z}_n > 0 \}} = \tilde{Z}_n {\bf 1}_{\{ \tilde{Z}_n > 0 \}}, \quad
\mbox{\rm where} \quad \tilde{Z}_n := \frac{V_n}{A_0 \cdots A_n}.
\end{equation}
\end{Lm}

\Pf
It follows by an inductive argument
that for all $n \in {\mathbb N}$,
\begin{equation}{\label{lm5.3.1}}
V_n =\max \left(\ \sum_{i=0}^n B_i \prod_{j=i+1}^n A_j, \bigvee_{k=1}^n \left[  \sum_{i=k}^n B_i \prod_{j=i+1}^n A_j + D_k \prod_{j=k}^n A_j
\right] \right),
\end{equation}
where $B_0 := V_0$.  Hence
\begin{equation} \label{NewLm5.4.1}
\tilde{Z}_n :=
\frac{V_n}{A_0 \cdots A_n} = \sum_{i=0}^n \frac{B_i}{A_0 \cdots A_i}
  - \min \left(0, {\mathfrak M}_n \right),
\end{equation}
where
\begin{equation} \label{NewLm5.4.1a}
{\mathfrak M}_n := \bigwedge_{k=1}^n \left[ \sum_{i=0}^{k-1}
      \frac{B_i}{A_0 \cdots A_i} - \frac{D_k}{A_0 \cdots A_{k-1}}\right].
\end{equation}

Now
\begin{equation} \label{NewLm5.4.2}
Z_n^{(p)} = F_{Y_0}^{(p)} \circ \cdots \circ F_{Y_n}^{(p)}(0) = \sum_{i=0}^n \frac{B_i}{A_0 \cdots A_i}.
\end{equation}
To obtain a similar expression for $Z_n^{(c)}$, observe
again by induction that
\begin{equation} \label{NewLm5.4.3}
Z_n^{(c)} = \min
\left(\sum_{i=0}^n \frac{B_i}{A_0 \cdots A_i},
\bigwedge_{k=0}^n
\left[ \sum_{i=0}^k \frac{B_i}{A_0 \cdots A_i} + \frac{D_k^\ast}{A_0
    \cdots A_k} \right] \right).
\end{equation}
Then substituting $D_0^\ast=-B_0$ and
$D_i^\ast = -A_iD_i - B_i$, $i=1,2,\ldots,$
yields that with ${\mathfrak M}_n$ given as in \eqref{NewLm5.4.1a},
\begin{equation} \label{NewLm5.4.4}
Z_n^{(c)} = \min \left(Z_n^{(p)}, 0 , {\mathfrak M}_n \right).
\end{equation} 
Comparing these expressions with \eqref{NewLm5.4.1}, we conclude that
\begin{equation} \label{NewLm5.4.6}
\tilde{Z}_n = Z_n^{(p)} - \min(0,{\mathfrak M}_n).
\end{equation}

Finally observe that if $\tilde{Z}_n \ge 0$, then
it follows from the previous equation that $0 \le
\max\left(Z_n^{(p)},Z_n^{(p)}-{\mathfrak M}_n\right)$
and hence
\[
\hspace*{1.8cm}
\tilde{Z_n} {\bf 1}_{\{ \tilde{Z}_n \ge 0 \}} =
\max\left(0,Z_n^{(p)},Z_n^{(p)}-{\mathfrak M}_n \right) = Z_n^{(p)} -
Z_n^{(c)}.
\hspace*{1.8cm}\halmos
\]

In the next lemma, we study the convergence of the sequence $\{
\tilde{Z}_n \}$ and relate it to the random variable $\bar{Z}^{(p)}$
defined just prior to the statement of Theorem 2.2.

%
%
%  LEMMA 5.5
%
%

\begin{Lm}
Assume the conditions of Lemma 5.2,
and let $\{ \tilde{Z}_n \}$ be defined as in \eqref{prelm3-1}.
Then in
$\mu_\xi$-measure, $\{ \tilde{Z}_n \}$ has the following regularity properties{\rm :}\\
\indent
{\rm (i)}  $\tilde{Z}_n \to \tilde{Z} $ a.s.\ as $n \to \infty$,
where $\tilde{Z}$ is a proper random variable supported on
$(0,\infty)$.
As a consequence, $Z^{(p)}$ and
$Z^{(c)}$ are well defined and $Z^{(p)} > Z^{(c)}$ a.s.\\
\indent
{\rm (ii)}   
${\bf E}_\xi \big[ \big( \bar{Z}^{(p)}\big)^\xi  \big] < \infty$.  
Moreover, for all $n$ and $u$,
\begin{equation} \label{lm-5.4.1S}
\left|\tilde{Z}_n {\bf 1}_{\{ n <  \tau \}} \right| \le \bar{Z}^{(p)}
   \quad \quad\mbox{\rm and}\quad 
\left|\tilde{Z}_{T_u} {\bf 1}_{\{ T_u <  \tau \}} \right| \le \bar{Z}^{(p)}.
\end{equation}
\end{Lm}

\Pf
We begin by establishing (ii).  Since
\begin{equation} \label{pf-lm5.4.1}
|V_n| \le A_n |V_{n-1}| + \left( A_n |D_n| + |B_n| \right),
\end{equation}
the process $\{ | V_n | \}$ is bounded from above by $\{ R_n \}$,
where $R_0 = |V_0|$ and, for each $n \in \pintegers$,
\begin{equation} \label{def-tildeB}
R_n = \tilde{B}_n + A_n R_{n-1}, \quad \mbox{\rm where} \quad \tilde{B}_n = \left( | B_n | + A_n |D_n| \right).
\end{equation}
Iterating the previous equation yields
\begin{equation} \label{pf-lm5.4.2}
R_n = \sum_{i= 0}^n \tilde{B}_i \prod_{i=j+1}^n A_i, \quad n=0,1,\ldots,
\end{equation}
where $\tilde{B}_0 \equiv |V_0 |$.
Since $| V_n | \le R_n$ for all $n$, we may now apply \eqref{prelm3-1}
to obtain that $\{ |\tilde{Z}_n| \}$ is bounded from above by the perpetuity sequence
\begin{equation} \label{pf-lm5.4.3}
W_n := (A_1 \cdots A_n)^{-1} R_n = \sum_{i=0}^n \frac{\tilde{B}_i}{A_1 \cdots A_i},
\end{equation}
and hence
\begin{equation}
|\tilde{Z}_n | {\bf 1}_{\{ n < \tau \}} \le W_n {\bf 1}_{\{ n < \tau
  \}} \le \sum_{i=0}^\infty  \frac{\tilde{B}_i}{A_0 \cdots A_i} {\bf
  1}_{\{ \tau > i \}} := \bar{Z}^{(p)},
\end{equation}
and similarly for $\tilde{Z}_{T_u} {\bf 1}_{\{ T_u < \tau \}}$.

It remains to show that ${\bf E} \big[ \big( \bar{Z}^{(p)} \big)^\xi
\big] < \infty$.  To this end, observe that if $\xi \ge 1$,
then by Minkowskii's inequality followed by a change of measure argument,
\begin{equation} \label{pf-lm5.4.4}
{\bf E}_\xi  \left[ \big( \bar{Z}^{(p)} \big)^\xi \right]^{1/\xi}
\le 
%\Bigg( {\bf E}_\xi  \Bigg[ \Bigg( \sum_{n=0}^\infty
%\left(A_1^{-1}Ê\cdots A_n^{-1}\right)\tilde{B}_n 
%{\bf 1}_{\{ \tau > n \}} % \Bigg)^\xi \Bigg] \Bigg)^{1/\xi}\\ \nonumber
%&  \le \sum_{n=0}^\infty {\bf E}_\xi \left[ \left(A_1^{-1} \cdots
%A_n^{-1}\right)^\xi 
%\tilde{B}_n^\xi {\bf 1}_{\{ \tau > n \}} \right]^{1/\xi}  \\ \nonumber
  {\bf E} \big[ \tilde{B}^\xi \big]^{1/\xi} \left( 1 +
    \sum_{i=1}^\infty {\bf P} \Big( \tau > i-1 \Big)^{1/\xi}
  \right) < \infty,
\end{equation}
where finiteness is obtained from $(H_2)$ and Lemma 5.1 (iv).
On the other hand,
if $\xi < 1$, then an analogous result is obtained using the
deterministic inequality \eqref{ex1-13} in place of Minkowskii's
inequality, which then completes the proof of (ii).

To establish (i), set 
$\tilde{Z}_- = \liminf_{n \to \infty} \tilde{Z}_n$ 
and observe that if ${\bf P} ( \tilde{Z}_- \le 0) = p > 0$,
then by Fatou's lemma 
\begin{equation} \label{limZ}
p \le \liminf_{n \to \infty} {\bf E} 
\left[ {\bf 1}_{\{ \tilde{Z}_n \le 0 \}} \right] = 
\liminf_{n \to \infty} {\bf P} \left( V_n \le 0 \right),
\end{equation}
where, in the last step, we have used Lemma 5.4 to replace $\tilde{Z}_n$ by $V_n$.
But by Lemma 5.2, $V_n \uparrow + \infty$ w.p.1 as $n \to \infty$, and thus \eqref{limZ}
is impossible.  We conclude $\tilde{Z}_- \in (0,\infty)$.

Next observe by definition of $\{ V_n \}$ that if $V_{n-1} \ge 0$,
\begin{equation} \label{ADDfeb16-1}
A_n V_{n-1} - |B_n| \le V_n \le A_n V_{n-1} + \left( |B_n| + A_n |D_n| \right).
\end{equation}
Let $\delta(n)$ denote the indicator function on the event $\{ V_l \ge 0,$ for all $l \ge n \}$.
Then by Lemma 5.2,
 $\lim_{n \to \infty}\delta(n) = 1$ a.s.
By iterating 
the left- and right-hand sides of \eqref{ADDfeb16-1}, we obtain that for any $m > n$,
\begin{equation} \label{pf-lm5.4.6}
\left| V_m - \left(A_{n+1} \cdots A_m \right) V_n \right| \delta(n) \le  \sum_{i= n+1}^m \tilde{B}_i \prod_{i=j+1}^n A_i,
\end{equation}
where $\tilde{B}_i$ is defined as in \eqref{def-tildeB}.
Hence for any $m > n$,
\begin{equation} \label{pf-lm5.4.7}
\left| \tilde{Z}_m - \tilde{Z}_n \right| \delta(n) \le \sum_{i=n+1}^\infty \frac{\tilde{B}_i}{A_1 \cdots A_i}.
\end{equation}

Now recall that  ${\bf E}_\xi \left[ \log A \right] = \Lambda^\prime(\xi) > 0$.  Hence the process
\[
\tilde{S}_n := -\log A_1 - \cdots -\log A_n + \log \tilde{B}_n
\]
has a negative drift, and using $(H_1)$ we have ${\bf E}_\xi \left[ (A^{-1})^\alpha \right] < \infty$ 
for some $\alpha > 0$.
Hence by Cram\'{e}r's large deviation theorem 
(\shortciteN{ADOZ98}, Section 2.2), it follows that
\begin{equation} \label{pf-lm5.4.8}
{\bf P}_\xi \left( \tilde{S}_n > n\left( -{\bf E}_\xi \left[ \log A \right] + \epsilon \right) \right) \le e^{-tn} \quad \mbox{\rm for all}\:\: n,
\end{equation}
for some constants $\epsilon \in \left(0,{\bf E}_\xi \left[ \log A \right] \right)$ and $t > 0$.  Thus
\begin{equation} \label{pf-lm5.4.9}
{\bf P}_\xi \left( \frac{ \tilde{B}_nÊ}{ A_1 \cdots A_n } > a^n \right) \le e^{-tn}, 
\end{equation}
where $a := \exp \left\{ -{\bf E}_\xi \left[ \log A \right] + \epsilon \right\} \in (0,1).$
Then by the Borel-Cantelli lemma,
\begin{equation} \label{pf-lm5.4.10}
{\bf P}_\xi \left( \frac{ \tilde{B}_nÊ}{ A_1 \cdots A_n } > a^n \:\: \mbox{\rm i.o.}\right) =0.
\end{equation}
Since $\lim_{n \to \infty} \delta(n)=1$ a.s.,
it follows as a consequence of \eqref{pf-lm5.4.7} and \eqref{pf-lm5.4.10} that $\tilde{Z}_n$ converges a.s.\ to a proper random
variable on $(0,\infty)$, which we denote by $\tilde{Z}$.

Since the above argument also holds for the recursion
$f(v)=Av+B$, it follows that $Z_n^{(p)}$ converges a.s.\ to a
random variable taking values in $(0,\infty)$.  Then by Lemma 5.4, the
limit in the definition of $Z^{(c)}$ must also exist and we must have
that
$Z^{(p)} - Z^{(c)} = \tilde{Z} > 0$ a.s.
\halmos

\subsection{Proofs of the main theorems}
The key to these proofs is the following result, which shows that $u^\xi{\bf E} \left[ N_u \right]$
behaves asymptotically as a product of two terms, the first being influenced mainly by the ``short-time" behavior of the process,
while the second term is determined by the ``large-time" behavior.  We will later identify this second term as a classical limit
related to the random walk process $S_n = \sum_{i=1}^n \log A_i$.

%
%
%  PROPOSITION 5.1
%
%

\begin{Prop}
Assume Letac's Model E, and suppose that $(H_1)$, $(H_2)$, and $(H_3)$ are satisfied.  Then
\begin{equation} \label{mainpf-1}
\lim_{u \to \infty} u^\xi {\bf E} \left[ N_u \right] = {\bf E}_\xi \left[ \tilde{Z}^\xi {\bf 1}_{\{ \tau=\infty \}} \right]
 \lim_{u \to \infty} {\bf E}_{\mathfrak D} \left[  N_u \left( \frac{V_{T_u}}{u} \right)^{-\xi} \Big| T_u < \tau \right].
\end{equation}
\end{Prop}

\Pf  Since $N_u  = g(V_0,\ldots, V_{\tau-1})$ for some 
function $g$, it follows by Lemma 5.3 that
\begin{equation} \label{pfprop1}
{\bf E} \left[ N_u \right] = {\bf E}_{\mathfrak D} \left[ N_u e^{-\xi S_{T_u}} \right].
\end{equation}
Moreover by the definition of $\tilde{Z}_n$ in \eqref{NewLm5.4.3},
\begin{eqnarray} \label{pfprop2}
e^{-\xi S_{T_u}} {\bf 1}_{\{T_u < \tau\}}&=& \left(  \left(A_1 \cdots A_{T_u} \right)^{-\xi} V_{T_u}^\xi  {\bf 1}_{\{T_u < \tau\}} \right) V_{T_u}^{-\xi} \\
\nonumber
  &=& \left( \tilde{Z}_{T_u}^\xi {\bf 1}_{\{T_u < \tau\}} \right) V_{T_u}^{-\xi}.
\end{eqnarray}
Combining the last two equations yields
\begin{equation} \label{pfprop3}
{\bf E} \left[ N_u \right] = {\bf E}_{\mathfrak D} \left[  \left(\Zt_{T_u}^\xi {\bf 1}_{\{ T_u < \tau\}} \right)   N_u  V_{T_u}^{-\xi} \right].
\end{equation}
Consequently, by conditioning on ${\mathfrak F}_{T_u \wedge (\tau-1)}$  and rearranging terms, we have that
\begin{equation} \label{pfprop4}
u^{\xi} {\bf E} \left[ N_u \right] = {\bf E}_{\mathfrak D} \left[ \left( \Zt_{T_u}^\xi {\bf 1}_{\{ T_u < \tau \}} \right) {\cal Q}_u \right],
\end{equation}
where 
\begin{equation} \label{pfprop4a}
{\cal Q}_u := {\bf E}_{\mathfrak D} \left[ N_u \left| {\mathfrak F}_{T_u \wedge (\tau-1)} \right. \right] \left( \frac{V_{T_u}}{u} \right)^{-\xi}
{\bf 1}_{\{ T_u < \tau \}}.
\end{equation}

Note
\begin{align} \label{pfprop5}
{\bf E}_{\mathfrak D} & \left[ \left( \Zt_{T_u}^\xi {\bf 1}_{\{ T_u < \tau \}} \right) {\cal Q}_u \right] 
 = {\bf E}_{\mathfrak D} \left[ \Zt_n^\xi {\cal Q}_u {\bf 1}_{\{ n \le T_u < \tau \}} \right] \\
\nonumber & \hspace*{3.5cm}
     + {\bf E}_{\mathfrak D} \left[ \left( \Zt_{T_u}^\xi {\bf 1}_{\{ T_u < \tau \}} - \Zt_n^\xi {\bf 1}_{\{ n \le T_u < \tau \}} \right) {\cal Q}_u  \right].
\end{align}
Now take the limit first as $u \to \infty$ and then as $n \to \infty$.

We begin by analyzing the first term on the right-hand side.
By Theorem 4.2 and by the definition of ${\cal Q}_u$ in \eqref{pfprop4a},
$\{ {\cal Q}_u \}$ is uniformly bounded in $u$.   Also, by Lemma 5.5 (ii),
${\bf E}_\xi \big[ |\Zt_n|^\xi {\bf 1}_{\{ n \le T_u  < \tau \}} \big] < \infty$ for all $n$.   Consequently for any given $n$,
\begin{equation} \label{pfprop6}
 \Ed \left[ \big|  \Zt_n^\xi {\cal Q}_u \big| {\bf 1}_{\{ n \le T_u < \tau \}} \right] \le K {\bf E}_\xi \left[ 
   \big|  \Zt_n \big|^\xi {\bf 1}_{\{ n \le T_u < \tau \}}  \right] \le K^\prime
\end{equation}
for finite constants $K$ and $K^\prime$.  In the middle term, we have replaced
$\Ed\left[ \cdot \right]$ with ${\bf E}_\xi \left[ \cdot \right]$, since these will be the same on the
set $\{ n \le T_u < \tau \}$.  Hence an application of the dominated
convergence theorem yields
\begin{align} \label{pfprop7}
&\lim_{u \to \infty} \Ed \left[ \Zt_n^\xi {\cal Q}_u {\bf 1}_{\{ n \le T_u < \tau \}} \right]
\\ \nonumber
& \hspace*{2cm} =  \Ed \left[ \Zt_n^\xi \lim_{u \to \infty}  {\bf 1}_{\{ n \le T_u < \tau \}}\Ed
 \left[ {\cal Q}_u \left| {\mathfrak F}_n,\, T_u < \tau \right. \right]  \right]
\\ \nonumber
& \hspace*{2cm}
= \Ed \left[ \Zt_n^\xi \lim_{u \to \infty} {\bf 1}_{\{ n \le T_u < \tau \}} \right] \Ed \left[ {\cal Q} \right],
\end{align}
where the last step of \eqref{pfprop7}
follows by Theorem 4.3 (and, in particular, the independence of the limiting distribution ${\cal Q}$ on the
initial distribution of the process $\{ V_k \}$, which is here taken
to be the distribution of $V_n$
for some fixed $n$).  In the second expectation on the right-hand side,
we have exchanged the limit and expectation, which is justified based on the weak convergence of ${\cal Q}_u$ to ${\cal Q}$ 
and the uniform boundedness of $\{ {\cal Q}_u \}$ (obtained from Theorem 4.2 and the definition of ${\cal Q}_u$).

To identify the first expectation on the right-hand side
of \eqref{pfprop7}, observe once again that 
$\Zt_n {\bf 1}_{\{ n \le T_u < \tau \}}$ is the same in the $\xi$-shifted measure as it is in the dual measure, since the dual measure
agrees with the $\xi$-shifted measure up until time $T_u$.  Hence
\[
\Ed \left[ \Zt_n^\xi \lim_{u \to \infty} {\bf 1}_{\{ n \le T_u < \tau \}} \right] = {\bf E}_\xi
  \left[ \Zt_n^\xi \lim_{u \to \infty} {\bf 1}_{\{ n \le T_u < \tau \}} \right] = {\bf E}_\xi \left[ \Zt_n^\xi {\bf 1}_{\{ \tau = \infty \}} \right],
\]
since an elementary argument yields that $T_u \uparrow \infty$ a.s.\ as $u \to \infty$.

Substituting into \eqref{pfprop7} and now taking the limit as 
$n \to \infty$ yields
\begin{align*}% \label{pfprop8}
\lim_{n \to \infty} \lim_{u \to \infty} \Ed \left[ \Zt_n^\xi {\cal Q}_u 
  {\bf 1}_{\{ n \le T_u <  \tau \}} \right] 
& = \lim_{n \to \infty} {\bf E}_\xi \left[ \Zt_n^\xi {\bf 1}_{\{ \tau = \infty \}} 
   \right] \Ed \left[ {\cal Q} \right]\\ 
& = {\bf E}_\xi \big[ \tilde{Z}^\xi {\bf 1}_{\{ \tau = 
  \infty\}} \big] \Ed \left[ {\cal Q} \right],
\end{align*}
where the last step follows from the dominated convergence theorem and
Lemma 5.5 (i) and (ii).
Finally, observe by Theorem 4.3 that
\begin{align*} 
\Ed \left[ {\cal Q} \right] & = \lim_{u \to \infty} \Ed \left[ {\cal Q}_u  \right] \left( {\bf P} \left( T_u < \tau \right) \right)^{-1}\\
& = \lim_{u \to \infty} \Ed \left[ 
{\bf E}_{\mathfrak D} \left[ N_u \left| {\mathfrak F}_{T_u \wedge (\tau-1)} \right. \right] \left( \frac{V_{T_u}}{u} \right)^{-\xi} {\bf 1}_{\{ T_u < \tau \}} \right]
\left( {\bf P} \left( T_u < \tau \right) \right)^{-1} \\ 
& = \lim_{u \to \infty} {\bf E}_{\mathfrak D} \left[  N_u \left( \frac{V_{T_u}}{u} \right)^{-\xi} \, \Big| \, T_u < \tau \right].
\end{align*}
[In the second equality, we have used that $N_u = 0$ on the set $\{ T_u \ge \tau \}$.]
Substituting the previous equation into \eqref{pfprop7} yields
\begin{align} \label{pfprop9}
\lim_{n \to \infty} & \lim_{u \to \infty}
  \Ed \left[ \Zt_n^\xi {\cal Q}_u 
{\bf 1}_{\{ n \le T_u \wedge (\tau-1) \}} \right]\\ \nonumber
&= {\bf E}_\xi \big[ \tilde{Z} {\bf 1}_{\{ \tau = \infty \}} \big]
 \lim_{u \to \infty} {\bf E}_{\mathfrak D}
 \left[  N_u \left( \frac{V_{T_u}}{u} \right)^{-\xi} \Big| T_u < \tau \right]. 
\end{align}
Going back to \eqref{pfprop4} and \eqref{pfprop5},
we see that the proof of the proposition will now be complete, provided that we can show that
\begin{equation} \label{pfprop10}
\lim_{n \to \infty} \lim_{u \to \infty}
{\bf E}_{\mathfrak D} \left[ \left( \Zt_{T_u}^\xi {\bf 1}_{\{ T_u < \tau \}} - \Zt_n^\xi {\bf 1}_{\{ n \le T_u < \tau \}} \right) {\cal Q}_u  \right] = 0.
\end{equation}

But recall once again that
$\{ {\cal Q}_u \}$ is uniformly bounded.  Thus, to establish \eqref{pfprop10}, it is sufficient to show
\begin{equation} \label{pfprop10a}
\lim_{n \to \infty} \lim_{u \to \infty}
{\bf E}_\xi \left[ \left( \Zt_{T_u}^\xi {\bf 1}_{\{ T_u < \tau \}} - \Zt_n^\xi {\bf 1}_{\{ n \le T_u < \tau \}} \right)^+  \right] = 0
\end{equation}
and
\begin{equation} \label{pfprop10b}
\hspace*{-0.01cm}
\lim_{n \to \infty} \lim_{u \to \infty}
{\bf E}_\xi \left[ \left( \Zt_{T_u}^\xi {\bf 1}_{\{ T_u < \tau \}} - \Zt_n^\xi {\bf 1}_{\{ n \le T_u < \tau \}} \right)^-  \right] = 0.
\end{equation}
Note that in these last
expectations, we have again replaced $\Ed\left[ \cdot \right]$ with ${\bf E}_\xi \left[ \cdot \right]$, since these expectations
involve random variables on $\{ T_u < \tau \}$, and on that set these expectations are actually the same.

To establish  \eqref{pfprop10a}, first apply Lemma 5.5 (ii).  Namely observe that
the integrand on the left-hand side is dominated by $2\big(
\bar{Z}^{(p)} \big)^\xi$, which 
is integrable.  Hence, applying the dominated convergence theorem twice, first with respect to the limit in $u$ and 
then with respect to the limit in $n$, we obtain
\begin{equation} \label{pfprop11}
{\bf E}_\xi \left[ \left( \lim_{u \to \infty} \Zt_{T_u}^\xi {\bf 1}_{\{ T_u < \tau \}} - \lim_{n \to \infty} 
  \Zt_n^\xi {\bf 1}_{\{ \tau = \infty \}} \right)^+  \right] = 0,
\end{equation}
where, in the last equality, we have used that $T_u \uparrow \infty$ as $u \to \infty$.  This establishes \eqref{pfprop10a}.
The proof of \eqref{pfprop10b} is analogous,  so we omit the ${\rm details.} \halmos$

Next we identify second term on the right-hand side of \eqref{mainpf-1} by relating it to the classical Cram\'{e}r-Lundberg constant.
First recall that $\tau^\ast$ is a typical return time of the random walk $S_n = \sum_{i=1}^n \log A_i$
to the origin, while $\tau$ is a typical regeneration for the process $\{ V_n \}$ under study.

%
%
%  LEMMA 5.6
%
%

\begin{Lm}
Assume Letac's Model E, and suppose that $(H_1)$, $(H_2)$, and $(H_3)$
are satisfied.  Then
\begin{equation} \label{Lmb-main}
\frac{1- {\bf E} \big[ e^{\xi S_{\tau^\ast}} \big] }{{\bf E} \left[ \tau^\ast \right]}
\lim_{u \to \infty} {\bf E}_{\mathfrak D} \left[ N_u \left( \frac{V_{T_u}}{u} \right)^{-\xi}  \Big| T_u < \tau  \right] =  C^\ast,
\end{equation}
where $C^\ast$ is the Cram\'{e}r-Lundberg constant defined in \eqref{sr203}.
\end{Lm}

\Pf  Define
\[
W_n = \left(\log A_n + W_{n-1} \right)^+, \mbox{ for }n =0,1,\ldots.
\]
Then $\{ W_n \}$ is a random walk reflected at the origin, and since ${\bf E} \left[ \log A\right] < 0$,
it is well known that this process is a recurrent Markov chain which converges to a random variable $W$ whose distribution is the stationary distribution.
Also, by \shortciteN{DI72}, Lemma 1,
\begin{equation} \label{Lmb-pf1}
\lim_{u \to \infty} u^{-\xi} {\bf P} \left( W > \log u \right) = C^\ast,
\end{equation}
where $W := \lim_{n \to \infty} W_n$.

Set $\tau^\ast = \inf\big\{n \in \pintegers: W_n = 0 \big\}$.  Since $\{ W_n \}$
has an atom at the origin, $\tau^\ast +1$ is a
regeneration time of the Markov chain $\{ W_n \}$.  Hence by the representation formula in Lemma 2.3,
\begin{equation} \label{Lmb-pf2}
{\bf P} \left(W >  \log u \right) = \frac{{\bf E} \big[ N_u^\ast \big]}{{\bf E} [\tau^\ast]},
\end{equation}
where
\[
N_u^\ast := \sum_{n=1}^{\tau^\ast} {\bf 1}_{(\log u, \infty)} (W_n).
\]

Let $\mu_A$ denote the marginal distribution of $\log A$, and set
\[
\mu_{A,\xi} (E) =  \int_E e^{\xi x}d\mu_A(x), \quad \forall E \in {\cal B}(\reals).
\]
Let $T_u^\ast = \inf \left\{ n \in \pintegers:  W_n > \log u \right\}$,
and define the dual measure
\begin{equation} \label{DUAL-MEAS-a}
{\mathfrak L} (\log A_n) = \left\{ \begin{array}{l@{\quad}l}
\mu_{A,\xi} & \mbox{for}\:\:  n=1,\ldots,T_u^\ast, \\
\mu_A & \mbox{for}\:\: n > T_u^\ast.
\end{array} \right.
\end{equation}
With a slight abuse of notation, let ${\bf E}_{\mathfrak D} \left[ \cdot \right]$ denote expectation with respect to this dual measure.
Then by a change of measure, analogous to Lemma 5.3, we obtain
\begin{align} \label{Lmb-pf30}
{\bf E}  \big[ N_u^\ast \big] &=  u^{-\xi} {\bf E}_{\mathfrak D} \left[ N_u^\ast
e^{-\xi(W_{T^\ast_u} - \log u)} \right] \\
\nonumber & =  u^{-\xi} {\bf E}_{\mathfrak D} \left[
{\bf E}_{\mathfrak D} \left[N_u^\ast \left| {\mathfrak F}_{T^\ast_u \wedge \tau^\ast} \right. \right]
  e^{-\xi(W_{T_u^\ast} - \log u)}  {\bf 1}_{\{ T_u^\ast \le \tau^\ast
    \}} \right],
\end{align}
which is equal to
\[
u^{-\xi} {\bf E}_{\mathfrak D} \left[
{\bf E}_{\mathfrak D} \left[N_u^\ast \left| {\mathfrak F}_{T^\ast_u \wedge \tau^\ast} \right. \right]
  e^{-\xi(W_{T_u^\ast} - \log u)} \Big| T_u^\ast \le \tau^\ast \right] {\bf P}_{\mathfrak D} \left( T_u^\ast \le \tau^\ast \right).
\]

To complete the proof, notice that the multiplicate random walk process
$\{ \exp(W_n) \}$ also satisfies the conditions of Theorem 4.3, and the
random variable inside the last expectation converges to the same
weak limit as the corresponding quantity for the process $\{ V_n \}$.
(Note that we write ``$T_u^\ast \le \tau^\ast$'' while we write ``$T_u < \tau$,''
since $\tau$ is a regeneration time for the original process, while
$\tau^\ast + 1$ is a regeneration time for the reflected random walk process.)
Thus by Theorem 4.3, $\left( J_1(u)/J_2(u) \right) \to 1$ as $u \to \infty$,
where
\begin{align*}
& \hspace*{.75cm}J_1(u) := {\bf E}_{\mathfrak D} \left[ {\bf E}_{\mathfrak D} \left[ N_u \left| {\mathfrak F}_{T_u} \right. \right] \left( \frac{V_{T_u}}{u} \right)^{-\xi} \Big|
 T_u < \tau \right];\\
&J_2(u) := {\bf E}_{\mathfrak D} \left[
{\bf E}_{\mathfrak D} \left[N_u^\ast \left| {\mathfrak F}_{T^\ast_u \wedge \tau^\ast} \right. \right]
  e^{-\xi(W_{T_u^\ast} - \log u)} \Big| T_u^\ast \le \tau^\ast \right].
\end{align*}
Consequently we conclude
\begin{align} \label{Lmb-pf100}
&\lim_{u \to \infty} {\bf E}_{\mathfrak D} \left[ N_u \left( \frac{V_{T_u}}{u} \right)^{-\xi}  \Big| T_u < \tau  \right]
= \lim_{u \to \infty}  \frac{u^\xi {\bf E} \left[ N_u^\ast \right]}{
  {\bf P}_{\mathfrak D} \left( T_u^\ast \le \tau^\ast \right)}\\ \nonumber
& \hspace*{3.5cm}
 = \lim_{u \to \infty} \frac{C^\ast {\bf E} \left[ \tau^\ast \right]}{ {\bf P}_{\mathfrak D} \left(
    T_u^\ast \le \tau^\ast \right)} \quad \mbox{\rm by \eqref{Lmb-pf1}
   and \eqref{Lmb-pf2}.}
\end{align}
Finally observe that $\lim_{u \to \infty}{\bf P}_{\mathfrak D} \left(
    T_u^\ast \le \tau^\ast \right) = {\bf P}_{\mathfrak D} \left( T_u^\ast = \infty \right)$,
and by a change of measure argument,
\begin{align} \label{Lmb-pf101}
{\bf P}_{\mathfrak D}  \left( \tau^\ast = \infty \right) &= 1 - {\bf P}_\xi \left( \tau^\ast < \infty \right) 
= 1 - \sum_{n=1}^\infty {\bf E}_\xi \big[ {\bf 1}_{\{ \tau^\ast = n \}} \big]
 \\ \nonumber
& \hspace*{-.75cm} = 1 - \sum_{n=1}^\infty {\bf E} \big[ e^{\xi S_n} {\bf 1}_{\{ \tau^\ast = n \}} \big] 
  = 1 - {\bf E} \big[ e^{\xi S_{\tau^\ast}} \big].
\end{align}
The required result then follows from \eqref{Lmb-pf100}
and \eqref{Lmb-pf101}.
\halmos\\[-.35cm]

%
%
%  PROOF OF THEOREM 2.1
%
%

{\sc Proof of Theorem 2.1.}
By Lemma 2.3,
${\bf P} \left( V > u \right) = {\bf E} \left[ N_u \right]/{\bf E}
  \left[ \tau \right]$, 
and hence by Proposition 5.1,
\begin{align} \label{pfLM2.2.2}
&\lim_{u \to \infty} u^\xi {\bf P} \left( V > u \right) 
 = {\bf E}_\xi \left[ \tilde{Z}^\xi {\bf 1}_{\{ \tau=\infty \}} \right] \\ \nonumber
 & \hspace*{3.5cm} \cdot \left( {\bf E}\left[ \tau \right] \right)^{-1}
 \lim_{u \to \infty} {\bf E}_{\mathfrak D} \left[  N_u \left( \frac{ V_{T_u} }{u} \right)^{-\xi} \Big| T_u < \tau \right].
\end{align}
By Lemmas 5.4 and 5.5 (i), $\tilde{Z} = Z^{(p)} - Z^{(c)}$ a.s.,
and by Lemma 5.6, the last limit on the right-hand side of
\eqref{pfLM2.2.2} may be identified as
\[
\hspace*{3.1cm}
\frac{C^\ast {\bf E} \left[ \tau^\ast \right]}{1 - {\bf E}
  \left[ e^{\xi S_{\tau^\ast}} \right]} = \frac{1}{\xi \lambda^\prime(\xi)}
\quad\quad {\rm by}\quad\!\! \eqref{sr203}.
\hspace*{3.1cm} \halmos
\]

%
%
%  PROOF OF THEOREM 2.2
%
%

{\sc Proof of Theorem 2.2.}
A repetition of \eqref{pfprop4} and \eqref{pfprop4a} yields
\begin{align} \label{pf-thm2.2.1}
u^\xi {\bf E} \left[ N_u \right] &= {\bf E} \left[ \left( \tilde{Z}_{T_u}^\xi {\bf 1}_{\{ T_u < \tau \}} \right) {\bf E} \left[ N_u \left| {\mathfrak F}_{T_u \wedge \tau-1} \right. \right]
\left( \frac{V_{T_u}}{u} \right)^{-\xi} \right] \\ \nonumber
& \le {\bf E} \left[ \big( \bar{Z}^{(p)} \big)^\xi \left( C_1(u) \log\left( \frac{V_{T_u}}{u} \right) + C_2(u) \right) \left( \frac{V_{T_u}}{u} \right)^{-\xi} \right],
\end{align}
where the last step follows from Theorem 4.2 and Lemma 5.5 (ii).
Moreover, since $V_{T_u} \ge u$, we obviously have
\begin{align} \label{primeC}
&\left( C_1(u) \log \left( \frac{V_{T_u}}{u} \right) + C_2(u) \right) \left( \frac{V_{T_u}}{u}\right)^{-\xi}\\ \nonumber
& \hspace*{3cm} \le \sup_{z \ge 0} \Big\{ e^{-\xi z}\left( zC_1(u) + C_2(u) \right)\! \Big\},
\end{align}
which is bounded.
Substituting this deterministic bound into the right-hand side of \eqref{pf-thm2.2.1} establishes the theorem.
The limiting values of $C_1(u)$ and $C_2(u)$ are obtained by a further
application of Theorem 4.2.
\halmos

Returning to the extremal index described in \eqref{ex3-6},
recall that
\begin{equation} \label{extremes1}
\Theta = \lim_{u \to \infty} \frac{{\bf P} \left( V_n > u,
\:\:\mbox{for some}\:\: n < \tau \right)}{ {\bf E} \left[N_u \right] }
= \lim_{u \to \infty} \frac{{\bf E} \left[ {\bf 1}_{\{ T_u < \tau \}} \right]}{
  {\bf E} \left[ N_u \right]}.
\end{equation}
Moreover one can show, similar to the proof of Proposition 5.1, that
\[
\lim_{u \to \infty} u^\xi {\bf E} \left[ {\bf 1}_{\{ T_u < \tau \}} \right] = {\bf E}_\xi \left[ \tilde{Z}^\xi {\bf 1}_{\{ \tau=\infty \}} \right]
 \lim_{u \to \infty} {\bf E}_{\mathfrak D} \left[  {\bf 1}_{\{ T_u < \tau \}} \left( \frac{V_{T_u}}{u} \right)^{-\xi} \Big| T_u < \tau \right].
\]
To identify the last term on the right-hand side, let
$T_u^\ast$ and $\tau^\ast$ represent the corresponding
terms in the reflected random walk model
$\{ W_n \}$ (as in the proof of Lemma 5.6).  Then by 
applying Lemma 4.1 to the process $\{ V_n \}$ and to the multiplicative
random walk $\{ \exp(W_n) \}$, we see that
$V_{T_u}/u$ and $\log W_{T_u^\ast}/\log u$ converge to the
same weak limit.  Hence
\begin{align}  \label{extremes2}
\lim_{u \to \infty} {\bf E}_{\mathfrak D} \left[\left. \left( \frac{V_{T_u}}{u}
  \right)^{-\xi} \right| T_u < \tau \right]
= \lim_{u \to \infty} {\bf E}_{\mathfrak D} \left[ \left. e^{-\xi(W_{T_u^\ast}
 - \log u)} \right| T_u^\ast \le \tau^\ast \right].
\end{align}

But the quantity on the right-hand side is just the Cram\'{e}r-Lundberg
constant $C^\ast$.  To see that this is the case, note that if $T_u^\ast
\not\le \tau^\ast$, then the unconstained random walk $S_n = \sum_{i=1}^n \log A_i$ returns to
the set $(-\infty,0]$, and starting from this new level, the overjump distribution
will have the same limiting distribution as it had when starting its first
cycle from the origin.  Thus, using the fact that the dual and $\xi$-shifted
measures are actually the same in this case, we obtain
\begin{align*} 
 \lim_{u \to \infty} {\bf E}_{\xi} \left[ \left. e^{-\xi(W_{T_u^\ast}
 - \log u)} \right| T_u^\ast \le \tau^\ast \right]
&= \lim_{u \to \infty} {\bf E}_{\xi} \left[ \left. e^{-\xi(S_{T_u^\ast}
 - \log u)} \right| T_u^\ast < \infty \right]\\
& = \lim_{u \to \infty}
u^\xi {\bf P} \left( T_u^\ast < \infty \right) = C^\ast,
\end{align*}
where the second equality follows from a change of measure argument
and the fact that ${\bf P}_\xi (T_u^\ast < \infty) = 1$.  Thus we conclude that
\begin{equation} \label{extremes301}
\lim_{u \to \infty} u^\xi {\bf E} \left[ {\bf 1}_{\{ T_u < \tau \}} \right]
  = C^\ast{\bf E}_\xi \left[ \tilde{Z}^\xi {\bf 1}_{\{ \tau=\infty \}} \right].
\end{equation}

Now by Proposition 5.1 and Lemma 5.6, we also have that
\begin{equation} \label{extremes4}
\lim_{u \to \infty} u^\xi {\bf E} \left[ N_u \right]
  = {\bf E}_\xi \left[ \tilde{Z}^\xi {\bf 1}_{\{ \tau=\infty \}} \right]
  \frac{C^\ast {\bf E} \left[\tau^\ast\right]}{1-{\bf E} \left[ e^{\xi
  S_{\tau^\ast}} \right]}.
\end{equation}
Substituting \eqref{extremes301} and \eqref{extremes4} into \eqref{extremes1},
we conclude \eqref{ex3-6}.

%
%
%  PROOFS OF THE RESULTS FROM NONLINEAR RENEWAL THEORY
%
%
\section{Proofs of the results from nonlinear renewal theory} 
 First we turn to the proof of Lemma 4.1, which extends a classical result
from nonlinear renewal theory (given in \shortciteN{MW82}, Theorem 4.2)
to the setting of our problem.\\[-.3cm]

%
%
%  PROOF OF LEMMA 4.1
%
%
{\sc Proof of Lemma 4.1.}
By definition of $\tilde{Z}_n$ in \eqref{prelm3-1},
\begin{equation} \label{pf-Lm4.1.1}
V_n = (A_1 \cdots A_n) \left( \tilde{Z}_n {\bf 1}_{\{ \tilde{Z}_n > 0 \}} + \tilde{Z}_n {\bf 1}_{\{ \tilde{Z}_n \le 0 \}} \right).
\end{equation}
Now by Lemma 5.2, $V_n \uparrow \infty$ w.p.1 under the measure $\mu_\xi$.  Thus $T_u < \infty$ a.s., and at this exceedance time, we obviously have 
${\bf 1}_{\{ \tilde{Z}_n \le 0 \}}=0$.
Thus $V_{T_u} = V^\prime_{T_u}$ for
\begin{equation} \label{pf-Lm4.1.2}
V_n^\prime = (A_1 \cdots A_n) Y_n, \quad n=1,2,\ldots,
\end{equation}
where $Y_n = Z_n$ on $\{ \tilde{Z}_n > 0 \}$ and $Y_n = 1$ otherwise,
and thus $Y_n$ is everywhere {\em positive}.
We begin by showing that
\begin{equation} \label{pf-Lm4.1.3}
\frac{V_{T_u}^\prime}{u} \Rightarrow \hat{V} \quad \mbox{\rm as} \quad u \to \infty,
\end{equation}
where $\hat{V}$ has the distribution described in \eqref{prelimNR2}.

To this end, note by \eqref{pf-Lm4.1.2} that
\begin{equation} \label{pf-Lm4.1.4}
\log V_n^\prime = S_n + \delta_n,
\end{equation}
where $S_n := \sum_{i=1}^n \log A_i$ and $\delta_n := \log Y_n$.  Hence $\{ \log V_n^\prime \}$ may be viewed
as a perturbed random walk where, under $\mu_\xi$-measure, $\{ S_n \}$ has a positive drift and the sequence 
$\{ \delta_n \}$ has the property that $\left\{ (\log A_i,\delta_i):  i=1,\ldots,n \right\}$ is independent of $\log A_j$
for all $j >n$.  This puts us in the setting of classical nonlinear renewal theory.

Next observe by Lemma 5.5 
and ${\bf 1}_{\{ \tilde{Z}_n \le 0 \}} \to 0$ a.s.
that $\delta_n:= \log Y_n$ converges a.s.\ to a proper random variable,
and hence $\delta_n/n
\to 0$ a.s.\ as $n \to \infty$.  Thus
$\{ \delta_n \}$ is slowly changing.  
By \shortciteN{MW82}, Theorem 4.2, it follows that
\eqref{pf-Lm4.1.3} holds, and hence this equation
 also holds with $V_{T_u}$ in place
of $V_{T_u}^\prime$.

To show that the result holds conditional on $\{ T_u < \tau \}$,
let $V_0 \sim \nu$, where $\nu$ is given in Lemma 2.2 (iv);
let $K_0=0$;  and let $K_1,K_2,\ldots$ denote the successive regeneration times.
For each $i$, let $R_i$ denote the distribution of
$(V_{T_u}/u)$ conditional on the event that $T_u \in [K_i,K_{i+1})$.
By independence of the regeneration cycles, 
$\{ R_i \}$ is an i.i.d.\ sequence of random variables.  Consequently, 
the conditional distribution of $V_{T_u}$ given that
$T_u \in [K_i,K_{i+1})$ is the same for all $i$.
Thus, the conditional distribution
of $V_{T_u}/u$ given $\{T_u < \tau \}$ must
agree with the unconditional distribution of $V_{T_u}/u$, and the
proof is complete. 
\halmos

%
%
%  PROOF:  PROPOSITION 6.1
%
%

Next we turn to the proofs of Theorems 4.1-4.3.
An important part of the proofs will be the study the process $\{ V_n \}$
as it returns from the rare set $(u,\infty).$  To do so,
we will introduce a new barrier at a lower level
$u^t$, where $t \in (0,1)$, and divide the trajectory into two parts,
namely that occuring prior to the event $V_n \in [-u^t,u^t]$,
and that occuring after this time.  
Intuitively, the process $\{ V_n \}$ will closely resemble a multiplicative random walk away from the set $[-u^{t},u^t]$,
and so a critical aspect of the proofs will be to characterize the behavior of $\{ V_n \}$ after it returns to $[-u^t,u^t]$,
but prior to regeneration.    Indeed, it is precisely in this region that $\{ V_n \}$ fails to resemble the random walk process.  In the next
proposition, we show that the behavior of the process 
in this critical region  may be neglected in an appropriate sense,
and we also provide a reasonably sharp estimate for the number of visits above the level $u$ which arise after returning to $[-u^t,u^t]$.

For any $v \ge 0$, first define
\[
K(v) = \inf\{ n:  |V_n | \le v \}.
\] 

\begin{Prop} Assume Letac's Model E, suppose that $(H_1)$,
$(H_2)$, and $(H_3)$ are satisfied, and assume that $V_0 > u$.
Then for any $t \in (0,1)$,
 there exist finite positive constants $\alpha$, $M$, and
 $\rho \in (0,1)$ such that 
\begin{equation} \label{newA3a}
{\bf E} \left[ \sum_{n=K(u^{t})}^{\tau -1} {\bf 1}_{\{ V_n > u\}}  \right] \le \Delta(u).
\end{equation}
The constant $\Delta(u)$ is characterized as follows.  Set
$\bar{V}_1 := A_1 \max \left(D_1,V_0 \right) + |B_1|$.  Then
\begin{equation} \label{newA3aa}
\Delta(u) := \frac{u^{-\alpha(1-t)}}{1-\rho} \left\{1+u^{-\alpha t}\!\! \sup_{w \in [-M,M]}\! {\bf E}_{w} \left[ \tau \right] {\bf E}_{M} \left[ \bar{V}_1^\alpha \right] \right\} < \infty
\end{equation}
{\rm (}where ${\bf E}_{w}\left[ \cdot \right]$ denotes 
expectation conditional on $V_0 = v${\rm )}.
\end{Prop}

Note that we have dropped the dependence on the dual measure in this proposition,
since we assume that $V_0 >u$, and hence the entire trajectory will take
place in the original measure.\\[-.3cm]

 {\sc Proof of Proposition 6.1.}  
 Recall by Lemma 5.1 (ii) that $\{ V_n \}$ satisfies a drift condition; namely, for some $\alpha > 0$,
 \begin{equation} \label{pf-lm4.3.1}
 {\bf E} \Big[ \left|V_n \right|^\alpha \big| V_{n-1}=v \Big] \le \rho |v|^\alpha
 + \beta {\bf 1}_{{\cal C}_M}(v), 
 \end{equation}
 where $\rho \in (0,1)$,
${\cal C}_M := [-M,M]$, and $M$ and $\beta$ are constants.  

We will divide the proof into three steps.  In the first step, we will
study the number of exceedances above level $u$ which occur
in an excursion beginning at time $K(u^t)$ ({\it i.e.}, when the process
$\{ |V_n| \}$ first falls below the level $u^t$ where $t<1$) and ending
at time $K(M)$ ({\it i.e.}, when this process first enters
the set ${\cal C}_M := [-M,M]$.  In the next step, we then
consider the number of exceedances above level $u$ which occur
between time $K(M)$ and the actual regeneration time.
Combining these two estimates in Step 3, we obtain the desired
upper bound.
\\[-.2cm]

%
%
%  STEP 1
%
%

{\it Step 1}:   For any $t \in (0,1)$,
\begin{equation} \label{pf-lm4.3.2}
{\bf E} \left[ \sum_{n=K(u^{t})}^{K(M)}
 {\bf 1}_{\{ V_n > u \}} \right] \le \frac{u^{- \alpha(1-t)}}{1-\rho}.
\end{equation}

{\it Proof}:  Suppose $V_0 = v$, where $|v| > M$.  Then iterating \eqref{pf-lm4.3.1} yields
\begin{equation} \label{pf-lm4.3.3}
{\bf E} \left[ |V_n|^\alpha {\bf 1}_{\{ K(M) > n\} } \big| V_0 = v \right] \le \rho^n |v|^\alpha, \quad n=0,1,\ldots,
\end{equation}
and hence
\begin{equation} \label{pf-lm4.3.4}
{\bf E} \left[ {\bf 1}_{\{ V_n > u \}} {\bf 1}_{\{ K(M) > n \}} \Big| V_0 = v \right] \le u^{-\alpha} \rho^n |v|^\alpha.
\end{equation}
Consequently,
\begin{equation} \label{pf-lm4.3.5}
{\bf E} \left[ \sum_{n=0}^{K(M)-1} {\bf 1}_{\{ V_n > u \}} \Big| V_0 = v \right] \le u^{-\alpha}  |v|^\alpha
\left(1 + \rho + \rho^2 + \cdots \right).
\end{equation}
Then by
the previous equation and the strong Markov property,
\begin{equation} \label{pf-lm4.3.6}
{\bf E} \left[ \sum_{n=K(u^t)}^{K(M)-1} {\bf 1}_{\{ V_n > u \}} \Big| V_{K(u^t)} \right] \le u^{-\alpha}  |V_{K(u^t)}|^\alpha
\left(1 + \rho + \rho^2 + \cdots \right),
\end{equation}
Since $|V_{K(u^t)}| \le u^t$ by definition, this
establishes \eqref{pf-lm4.3.2}.\\[-.2cm]

%
%
%  STEP 2
%
%

{\it Step 2}:  We have
\begin{equation}
{\bf E} \left[ \sum_{n=K(M)}^{\tau-1}  {\bf 1}_{\{V_n> u\}}
  \right]
 \le \frac{u^{-\alpha}}{1-\rho} \left\{ \sup_{w \in [-M,M]}
  {\bf E}_w \left[ \tau \right] {\bf E}_M \left[ \bar{V}_1^\alpha \right] \right\}.
\end{equation} 

{\it Proof}: 
It suffices to show that for any $v \in [-M,M]$,
\begin{equation} \label{pf-lm4.3.7}
{\bf E} \left[ N_u \big| V_0 = v \right] 
 \le \frac{u^{-\alpha}}{1-\rho} \left\{ \sup_{w \in [-M,M]} {\bf E}_w \left[ \tau \right] {\bf E}_M \left[ \bar{V}_1^\alpha \right] \right\}.
\end{equation}

To this end, introduce the augmented chain described in
Remark 2.1.   Namely consider the process $\{ (V_n,  \eta_n) \}$,
where $\{ \eta_n \}$ is an i.i.d.\ sequence of Bernoulli random variables with ${\bf P} \left( \eta_n =1 \right) = \delta$ and ${\bf P} \left( \eta_n=0 \right) = 1-\delta$.
As discussed in Remark 2.1, we may then identify $\tau$ as the return time of this augmented chain to the set ${\mathfrak C} := {\cal C} \times \{1\}$, where 
${\cal C}$ and $\delta$ are obtained from the minorization condition $({\cal M})$ (which holds by Lemma 5.1 (iii)). 
We assume without loss of generality that ${\cal C} \subset {\cal C}_M$,
and define ${\mathfrak M} = {\cal C}_M \times \{ 0,1 \}$.

With a slight abuse of notation, let
$P$ denote the transition kernel of the augmented chain $\{ (V_n,\eta_n) \}$.
Then introduce the taboo transition kernel 
\[
{}_GP(x,E) := \int_E {\bf 1}_{G^c}(y) P(x,dy),
\]
where $G$ and $E$ are Borel subsets of
$\reals \times \{ 0,1 \}$.
Then by the last-exit decomposition ({\it cf.} \shortciteN{SMRT93}, 
Section 8.2),
\begin{equation} \label{lastexit1}
{}_{{\mathfrak C}} P^n({\bf v},E) = {}_{{\mathfrak M}} P^n({\bf v},E)
 + \sum_{k=1}^{n-1} \int_{{\mathfrak M}-{\mathfrak C}} \: {}_{{\mathfrak C}}
P^k({\bf v},d{\bf w}) _{\mathfrak M} \! P^{n-k}({\bf w},E).
\end{equation}

Now set $E = (u,\infty) \times \{0,1\}$ and fix ${\bf v}=(v,q)$,
where $v < u$ and $q \in \{0,1\}$.  Then sum \eqref{lastexit1}
over all  $n \in \pintegers$.
On the left-hand side of \eqref{lastexit1}, the term 
${}_{\mathfrak C}P^n\big({\bf v},(u,\infty)\times\{0,1\}\big)$ describes the probability that regeneration does {\em not}\, occur
during the first $n$ time increments and that $V_n \in (u,\infty)$.  Thus 
\begin{align} \label{lastexit2}
&\sum_{n=1}^\infty {}_{\mathfrak C}P^n\big({\bf v},(u,\infty)\times\{0,1\}\big)\\ \nonumber
& \hspace*{1cm} = \sum_{n=1}^\infty {\bf E} \left[ {\bf 1}_{\{ V_n > u \}}
  {\bf 1}_{\{ \tau > n \}} \big| V_0 = v \right]
= {\bf E} \left[ N_u \big| V_0 = v \right].
\end{align}

Next consider the second term on the right-hand side of \eqref{lastexit1}.   
Applying Tonelli's theorem
to interchange the order of summation and integration, 
and then interchanging the order of summation, we obtain
\begin{align} \label{lastexit3}
\sum_{n=1}^\infty & \sum_{k=1}^{n-1}  \int_{{\mathfrak M}-{\mathfrak C}}\: {}_{\mathfrak C}P^k({\bf v},d{\bf w}) {}_{\mathfrak M} P^{n-k}\big({\bf w},(u,\infty)
\times \{0,1\} \big)\\
&=Ê\int_{{\mathfrak M}-{\mathfrak C}} \: \sum_{k=1}^\infty {}_{\mathfrak C}P^k({\bf v},
d{\bf w}) \left( \sum_{n=k+1}^\infty \: {}_{\mathfrak M} P^{n-k}\big({\bf w},
(u,\infty) \times \{0,1\} \big) \right) \nonumber\\
\nonumber
&=\int_{{\bf w} \in ( {\mathfrak M}-{\mathfrak C})} \int_{{\bf z} \in {\mathfrak M}^c}
\sum_{k=1}^\infty\: {}_{\mathfrak C}P^k  ({\bf v},d{\bf w})
{}_{\mathfrak M} P({\bf w},d{\bf z})\\[-.2cm] \nonumber
& \hspace*{4.25cm} \cdot \left( \sum_{n=0}^\infty \;\; {}_{\mathfrak M}P^{n}\big(
{\bf z},(u,\infty) \times \{0,1\}\big) \right). 
\end{align}
For the last term in parentheses, note that  ${}_{\mathfrak M}P^{n}\big(
{\bf z},(u,\infty) \times \{0,1\}\big)$ describes the probability that
$\{ V_n \}$ avoids the set $[-M,M]$ during the first $n$ time increments
and that $V_n \in (u,\infty)$.  
Hence setting ${\bf z} = (z,r)$ yields, by \eqref{pf-lm4.3.5},
\begin{align} \label{lastexit4}
 \sum_{n=0}^\infty  \: {}_{\mathfrak M}P^{n}\big({\bf z},(u,\infty) \times \{0,1\}\big)& = {\bf E} \left[ \sum_{n=0}^{\infty} {\bf 1}_{\{ V_n > u \}}
   {\bf 1}_{\{ K(M) > n \}} \Big| V_0 = z \right] \\ \nonumber
& \hspace*{-1.5cm}  = {\bf E} \left[ \sum_{n=0}^{K(M)-1} {\bf 1}_{\{ V_n > u \}}  \bigg| V_0 = z  \right]  \le \frac{u^{-\alpha}|z|^\alpha}{1-\rho}.
\end{align}
Substituting this inequality into the previous equation, 
we see that the right-hand side of \eqref{lastexit3} is bounded above by
\begin{equation} \label{lastexit5}
\frac{u^{-\alpha}}{1-\rho} \sum_{k=1}^\infty \int_{{\mathfrak M}-{\mathfrak C}}
 {}_{\mathfrak C}P^k  ({\bf v},d{\bf w}) \left(  \int_{{\mathfrak M}^c} 
{}_{\mathfrak M} P({\bf w},d{\bf z}) |z|^{\alpha} \right).
\end{equation}
Note that with ${\bf w}=(w,s)$,
\[
 \int_{{\mathfrak M}^c} 
{}_{\mathfrak M} P({\bf w},d{\bf z}) |z|^{\alpha}
\le {\bf E} \left[ \bar{V}_1^\alpha \big| V_0 = M \right],
\]
where $\bar{V}_1 := A_1 \max \left(D_1,V_0 \right) + |B_1|$.
Moreover, for the remaining integral in \eqref{lastexit5}, we have
by definition that
\[
 \int_{{\mathfrak M}-{\mathfrak C}}
 {}_{\mathfrak C}P^k  ({\bf v},d{\bf w}) = 
{\bf P} \left( \tau > k, V_k \in [-M,M] \big| V_0 = v \right).
\]
Substituting the last two estimates into \eqref{lastexit5},
we obtain that \eqref{lastexit5} is bounded above by
\begin{equation} \label{lastexit6}
\frac{u^{-\alpha}}{1-\rho} \sum_{k=1}^\infty {\bf P} \left( \tau > k \, \big| V_0 = v \right) {\bf E} \left[ \bar{V}_{1}^\alpha \, \big| V_0 = M \right].
\end{equation} 

Now repeat the same argument, but applied to the first term on the right of \eqref{lastexit1}.  Essentially, this can be viewed
as one of the terms in the previous sum,
namely the term $k=0$.
In fact, 
the previous argument may be repeated without change to obtain
\begin{equation} \label{lastexit7}
{}_{{\mathfrak M}} P^n({\bf v},E) \le {\bf P} \left( \tau > 0 \, \big| V_0 = v \right) {\bf E} \left[ \bar{V}_{1}^\alpha \, \big| V_0 = M \right].
\end{equation}
Substituting these last two equations into the right-hand side of \eqref{lastexit1} and substituting \eqref{lastexit2}
into the left-hand side of \eqref{lastexit1} yields
\begin{align} \label{lastexit8}
{\bf E} \left[ N_u \big| V_0 = v \right] & = \frac{u^{-\alpha}}{1-\rho} \sum_{k=0}^\infty {\bf P} \left( \tau > k \, \big| V_0 = v \right) {\bf E} 
\left[ \bar{V}_{1}^\alpha \, \big| V_0 = M \right] \\ \nonumber
& \le \frac{u^{-\alpha}}{1-\rho} \left\{ \sup_{w \in [-M,M]} {\bf E}_w \left[ \tau \right] {\bf E}_M \left[ \bar{V}_1^\alpha \right] \right\},
\end{align}
which is \eqref{pf-lm4.3.7}.

Finally observe that the quantity on the right-hand side of
\eqref{lastexit8} is finite.  To this end,
it is sufficient to verify that
the set ${\mathfrak M}$ is regular.  For this purpose, observe
that since $\{ \eta_n \}$ is an i.i.d.\ sequence, it follows 
by a slight modification of Lemma
5.1 (iii) that ${\mathfrak M}$ is petite. 
Moreover, letting $\kappa$ denote the first return time of 
$\{ V_n \}$ to ${\cal C}_M$, then by Lemma 5.1 (i) and
Theorem 15.0.1 of \shortciteN{SMRT93}, we have that
$\sup_{v \in {\cal C}_M} {\bf E}_v \big[t^{\kappa} \big] < \infty$ for some
$t > 1$.
Then by definition, it follows that this 
last equation also holds for the augmented chain with ${\mathfrak M}
:= {\cal C}_M \times \{0,1\}$ in place
of ${\cal C}_M$.  Consequently, the conditions of \shortciteN{SMRT93},
Theorem 11.3.14 (i), are fulfilled and hence ${\mathfrak M}$ is regular.\\[-.2cm]

%
%
%  STEP 3
%
%

{\it Step 3}: 
Finally observe that by summing the expectations studied in Steps 1 and 2,
we immediately obtain \eqref{newA3a}.
\halmos\\[-.3cm]

%
%
%  PROOF:  THEOREM 4.1
%
%
{\sc Proof of Theorem 4.1.}
By Proposition 6.1, it is sufficient to show that for some $t \in (0,1)$,
\begin{equation} \label{pf-thm4.2.1}
{\bf E} \left[ \sum_{n=0}^{K(u^t)-1} {\bf 1}_{\{ V_n>u \}} \bigg| \frac{V_0}{u} = v \right] = U(\log v),
\end{equation}
where $U(z):= \sum_{n \in \integers} \mu_A^{\ast n}(-\infty,z)$ and $\mu_A$ is the marginal distribution
of $-\log A$.   As explained in Section 4, 
$U$ may be viewed under a minor continuity condition as the renewal function of $-S_n = -\sum_{i=1}^n \log A_i$,
while the expectation on the left-hand side may be viewed as a {\em truncated}\, renewal function of the nonlinear
process $\{ V_n \}$.

In the remainder of the proof, we will suppress the conditioning in \eqref{pf-thm4.2.1}.

To prove \eqref{pf-thm4.2.1}, we first establish an upper bound and then a corresponding lower bound.
For the upper bound, begin by observing ({\it cf.}\ \eqref{pf-drift1}) that
\begin{equation} \label{pf-thm4.2.2}
\frac{|V_n|}{|V_{n-1}|} \le A_n + \frac{\big(A_n |D_n| + |B_n| \big)}{|V_{n-1}|}.
\end{equation}
Note by definition that $|V_k| > u^t$ for all $k \le K(u^t)$.
Hence
\begin{equation} \label{pf-thm4.2.3}
\log \left( \frac{|V_n|}{|V_{n-1}|} \right) \le \log 
\bigg( A_n  + u^{-t} \Big(A_n |D_n| + |B_n| \Big)\bigg), \quad \mbox{\rm all} \:\: n < K(u^t).
\end{equation}
Now introduce the random walk
\begin{equation} \label{defSu}
S_n^{(u)} := \sum_{i=1}^n X_i^{(u)}, \:\:\mbox{\rm where}\:\: X_i^{(u)} :=  \log \left( A_i + u^{-t} \big(A_n |D_n| + |B_n| \big) \right),
\end{equation}
and $S_0^{(u)}=0$.
It follows from the previous two equations that
$\log |V_n| - \log |V_0| \le S_n^{(u)}$.  Since 
$(V_0/u) = v$, it then follows that
\begin{equation} \label{pf-zz1}
\log |V_n| - \log u \le S_n^{(u)} + \log v.
\end{equation}
Now suppose that $|V_n|>u$.  Then the left-hand side of \eqref{pf-zz1}
is positive, and thus $S_n^{(u)} > -\log v$.  Consequently,
\begin{equation} \label{pf4.2.5}
{\bf E} \left[ \sum_{n=0}^{K(u^t)-1} {\bf 1}_{\{ |V_n|>u \}} \right] \le {\bf E} \left[ \sum_{n=0}^\infty {\bf 1}_{\{S_n^{(u)} > -\log v \}} \right]
  := U^{(u)}(\log v).
\end{equation}

To establish an upper bound, it remains to show that
\begin{equation} \label{pf4.2.6}
\limsup_{u \to \infty} U^{(u)}(\log v) \le U(\log v).
\end{equation}
To this end, observe
by $(H_2)$ that ${\bf E} \big[ \exp \big\{ \epsilon X_i^{(u)} \big\} \big] < \infty$ for some $\epsilon > 0$ (since we take the logarithm
of the random variable $A_n |D_n| + |B_n|$ in the definition of $X_i^{(u)}$).  
Moreover,   
${\bf E} \left[ \log A \right] < 0 \Longrightarrow {\bf E} \big[ X_i^{(u)} 
\big] < 0$ for sufficiently large $u$.  
Since $U^{(u)}$ can itself be viewed as a renewal function
for the random walk $-S_n^{(u)}$, it follows that
$U^{(u)}(\log v) < \infty$.
Since $X_i^{(u)}$ decreases monotonically to $\log A_i$ as $u \to \infty$,
it follows by a dominated convergence argument that
\begin{equation} \label{pf4.2.7}
\lim_{u \to \infty} U^{(u)}(\log v) 
= {\bf E} \left[ \sum_{n=0}^\infty {\bf 1}_{\{ S_n > -\log v \}} \right] = U(\log v),
\end{equation}
which establishes the required upper bound.

Turning now to the lower bound, fix $\epsilon > 0$ and choose $N$ sufficiently large such that
\begin{equation} \label{pf4.2.8}
\left| {\bf E} \left[ \sum_{n=0}^N {\bf 1}_{\{ S_n > -\log v \}} \right]
  - U(\log v) \right| < \epsilon.
\end{equation}
Then
\begin{equation} \label{pf4.2.9}
{\bf E} \left[ \sum_{n=0}^{K(u^t)-1} {\bf 1}_{\{ V_n>u \}} \right] \ge {\bf E} \left[ \sum_{n=0}^{N} {\bf 1}_{\{ V_n>u \}} \right] - N{\bf P}
  \left( K(u^t) \le N \right).
\end{equation}
To bound the first term on the right-hand side, note that
\[
V_n \ge A_n V_{n-1} - |B_n|,
\]
and by iterating this equation we obtain
\begin{equation} \label{pf4.2.10}
V_n \ge (A_1 \cdots A_n)V_0 - \sum_{i=1}^{n-1} \prod_{j=i+1}^n A_j |B_i|, \quad \mbox{\rm for all}\:\: n.
\end{equation}
Thus, setting $(V_0/u) = v$, we obtain for any given $n$ that
\begin{equation} \label{pf4.2.11}
\liminf_{u \to \infty} \frac{V_n}{u} \ge \left( A_1 \cdots A_n \right) v := \exp \left\{ S_n + \log v \right\} \quad \mbox{\rm a.s.}
\end{equation}
Hence, if
$S_n > -\log v$, then we will necessarily have
$\liminf_{u \to \infty} (V_n/u)> 1$.  Consequently,
\begin{equation} \label{pf4.2.12}
\liminf_{u \to \infty} {\bf E} \left[ \sum_{i=1}^N {\bf 1}_{\{ V_n > u \}} \right] \ge {\bf E} \left[ \sum_{n=1}^N {\bf 1}_{\{ S_n > -\log v \}} \right]
\ge U(\log v)-\epsilon,
\end{equation}
where the last step was obtained by \eqref{pf4.2.8}.   Finally, observe by the definition of $K(u^t)$ and \eqref{pf4.2.10} that for
any $s \in (t,1)$,
\begin{align} \label{pf4.2.15}
{\bf P} \left( K(u^t) \le N \right) & \le {\bf P} \Big( \left(A_1 \cdots A_N\right)V_0 \le u^s \Big)\\ \nonumber
& \quad + {\bf P} \left(  \sum_{i=1}^{N-1} \prod_{j=i+1}^N A_j |B_i| > u^s - u^t \right).
\end{align}
In the first term on the right-hand side, $(V_0/u)= v$, and so the probability in question reduces to
${\bf P} \left( \left(A_1 \cdots A_N \right) v \le u^{s-1} \right) = {\bf P} \left( S_N + \log v \le (s-1)\log u \right)$, and the latter probability tends to zero
as $u \to \infty$,
as $N$ is fixed and $(s-1) < 0$.  Similar reasoning shows that the second term on the right of \eqref{pf4.2.15}  also tends to zero as
$u \to \infty$.  Thus we conclude $\lim_{u \to \infty} {\bf P} \left( K(u^t) \le N \right) = 0$ for any fixed $N$.
Substituting this last equation and \eqref{pf4.2.12} into \eqref{pf4.2.9}
yields the desired lower bound.
\halmos

%
%
%  PROOF OF THEOREM 4.2
%
%
 {\sc Proof of Theorem 4.2.}
 First assume that $(V_0/u) = v$ for some $v > 1$.  
Then it follows by \eqref{pf4.2.5}
 and Proposition 6.1 that
 \begin{equation} \label{pf-UB1}
 {\bf E} \left[ N_u  \left| \frac{V_0}{u} = v \right. \right] =
 {\bf E} \left[ \sum_{n=0}^{\infty} {\bf 1}_{\{ V_n>u \}} \left| \frac{V_0}{u} = v \right. \right] \le U^{(u)}(\log v) + \Delta(u),
 \end{equation}
 where
\begin{equation} \label{pf-UB2}
 U^{(u)} (\log v) :=  {\bf E} \left[ \sum_{n=0}^\infty {\bf 1}_{\{S_n^{(u)} > -\log v \}} \right]
 \end{equation}
 and $S_n^{(u)}$ is defined as in \eqref{defSu}.  In the definition of $U^{(u)}$, note that the terms inside the sum are obviously
 bounded above by ${\bf 1}_{\{ -S_n^{(u)} \le \log v \}}$, and substituting the latter quantity into the right-hand side of \eqref{pf-UB2},
 we obtain the renewal function of the process $\{ -S_n^{(u)} \}$.  
 Consequently, Lorden's inequality (\shortciteN{SA03}, Proposition V.6.2) yields
 \begin{equation} \label{pf-UB3}
 U^{(u)}(\log v) \le \frac{\log v}{m_u} + \left( 1+ \frac{\sigma_u^2}{m_u^2} \right),
 \end{equation}
 where
 \begin{equation} \label{pf-UB4}
 m_u := {\bf E} \big[ X^{(u)} \big] \quad \mbox{\rm and} \quad \sigma^2_u := {\rm Var}\: \big( X^{(u)} \big)
 \end{equation}
and $X^{(u)}$ is defined as in \eqref{defSu}.

Observe that these moments are actually finite.  In particular,
$|x|^\epsilon$ dominates $(\log x)^2$ for any $\epsilon > 0$.
Hence, applying $(H_2)$ for a sufficiently small choice of $\epsilon >0$,
we obtain
\begin{equation} \label{pf-UB4a}
{\bf E} \left[ \big( X^{(u)}\big)^2 \right] \le {\rm const.}\cdot {\bf E} \left[ e^{\epsilon X^{(u)}} \right] < \infty, \quad \mbox{\rm for all}\:\: u.
\end{equation}
Consequently, both the constants $m_u$ and $\sigma^2_u$ are finite.
Substituting \eqref{pf-UB3} into \eqref{pf-UB1} yields
\begin{equation} \label{pf-UB5}
 {\bf E} \left[ N_u  \left| \frac{V_0}{u} = v \right. \right] \le C_1(u) \log v + C_2(u),
 \end{equation}
 where ({\it cf.}\ \eqref{newA3aa}, \eqref{pf-UB4})
 \begin{equation} \label{pf-UB6}
 C_1(u):= \frac{1}{m_u} \quad \mbox{\rm and} \quad C_2(u) :=  \left(1+ \frac{\sigma_u^2}{m_u^2} \right) + \Delta(u).
 \end{equation}
 
 Next recall that our primary objective is to study $\Ed \left[ N_u \left| {\mathfrak F}_{T_u \wedge (\tau-1)} \right. \right].$
But by the strong Markov property,
\[
\Ed \left[ N_u \left| {\mathfrak F}_{T_u \wedge (\tau-1)} \right. \right] = \Ed 
\left[ N_u \big| V_{T_u} \right] {\bf 1}_{\{ T_u < \tau \}}.
\] 
Observe that in the dual measure, the process $\{ V_n \}$ reverts
back to its original measure after time $T_u$, and if $T_u \not< \tau$
then $N_u$ is zero.
Thus, we may apply \eqref{pf-UB5} to the right-hand side of the previous equation to obtain that
\begin{equation} \label{pf-UB7}
\Ed \left[ N_u \left| {\mathfrak F}_{T_u \wedge (\tau-1)} \right. \right] \le 
\left(C_1(u) \log \left(\frac{V_{T_u}}{u} \right) + C_2(u) \right) 
{\bf 1}_{\{ T_u < \tau \}},
\end{equation}
which is the required upper bound.

Finally observe that as $u \to \infty$, $X^{(u)}$ decreases monotonically to $\log A$ a.s.\ and $\Delta(u) \downarrow 0$, and hence
as $u \to \infty$,
\begin{equation} \label{pf-UB8}
C_1(u) \to C_1 := \frac{1}{m} \quad \mbox{\rm and} \quad C_2(u) \to C_2 :=1+\frac{\sigma^2}{m^2}
\end{equation}
where $m:= {\bf E} \left[ \log A \right]$ and $\sigma^2 := {\rm Var} \: (\log A)$.
\halmos

%
%
%  PROOF:  THEOREM 4.3
%
%
{\sc Proof of Theorem 4.3.}
By the strong Markov property,
\[
{\bf E} \left[ N_u \, \left| \, \frac{V_{T_u}}{u} = v \right. \right] = {\bf E} \left[ N_u \, \left| \, \frac{V_0}{u} = v \right. \right] := H_u(v),
\]
and by Lemma 4.1,
\begin{equation} \label{pf-Lm4.2.1}
\lim_{u \to \infty} H_u(v) = H(v):= U(\log v).
\end{equation}
Let $\nu_u$ denote the probability law of $V_{T_u}/u$ in the dual measure.
Then our objective is to show that for any $z > 1$,
\begin{equation} \label{pf-Lm4.2.2}
\lim_{u \to \infty} \int_1^z H_u(v) v^{-\xi} d\nu_u(v) = \int_1^z H(v) v^{-\xi} d\hat{\mu}(v).
\end{equation}
[Here we have taken the lower end point of the integral at one, since $V_{T_u}/u > 1$ for every $u$.]

Next observe that
\begin{equation}
\int_1^z H_u(v) v^{-\xi} d\nu_u(v) = R_1(u) + R_2(u),
\end{equation}
where
\[
R_1(u) = \int_1^z H_u(v) v^{-\xi} d\left(\nu_u(v) - \hat{\mu}(v) \right),\quad
 R_2(u) = \int_1^z H_u(v) v^{-\xi} d\hat{\mu}(v).
\]
Now by
\eqref{pf-UB5}, \eqref{pf-UB8}, 
\eqref{pf-Lm4.2.1},
and the dominated convergence theorem,
it follows that as $u \to \infty$,
\begin{equation}
R_2(u) \to R_2 := \int_1^z H(v) v^{-\xi} d\hat{\mu}(v).
\end{equation}
Finally, $R_1(u) \to 0$ as $u \to \infty$ using
\eqref{pf-UB5}, \eqref{pf-UB8}, and the weak
convergence of $\nu_u$ to
$\hat{\mu}$.
\halmos

\section{A generalization of Letac's Model E}
While the results of the previous sections were obtained specifically for Letac's Model E, the methods are 
actually quite general, applying to rather arbitrary recursions subject to appropriate regularity conditions.  However, if one departs from the
setting of the previous sections, then
one only obtains a recursive algorithm for computing the constant $C$, rather than the more elegant expression in terms of Letac's backward
recursive equation.

As an illustration of our method in a more general setting, consider the polynomial recursion
\begin{equation} \label{poly-1}
V \stackrel{\cal D}{=} F_Y(V) := B_k V + B_{k-1} V^{(k-1)/k} + \cdots + B_{1} V^{1/k} + B_0,
\end{equation}
where $(B_0,\ldots,B_k)$ is a sequence of positive-valued random variables and $Y=(B_0,\ldots,B_{k-1})$.  Note
\eqref{poly-1} may be written in the more suggestive form 
\begin{equation} \label{poly-2}
V \stackrel{\cal D}{=} A V + G_Y(V), \:\: \mbox{\rm where}\:\: A := B_k \:\: \mbox{\rm and} \:\: G_Y(V) := \sum_{j=0}^{k-1} B_j V^{j/k}.
\end{equation}
Observe that when the coefficients $(B_0,\ldots, B_{k-1})$ are fixed, $v^{-1}G_Y(v) \downarrow 0$ as $ v \to \infty$.  Thus, for large $V$,
the random function $G_Y(V)$ may be viewed as a remainder term, whereas the asymptotic behavior of the process is dominated 
in this case by the first term, $AV$.

Starting with the SFPE \eqref{poly-1},
we can form the forward recursion
\begin{equation} \label{poly-3}
V_n = A_n V_{n-1} + G_{Y_n}(V_{n-1}), \quad n=1,2,\ldots, \quad V_0=v,
\end{equation}
where $\{ A_n \}$ is an i.i.d.\ sequence of random variables having the same probability law as $A$,
and $\{ G_{Y_n} \}$ is an independent sequence
of random functions having the same law as $G_Y$ and
whose $n$th member will depend on $A_n$.  
Essentially, \eqref{poly-3} can be viewed as an extension of Example 3.2 of Section 3, where the term ``$B$" in the SFPE \eqref{ex2-3}
is now replaced with the random function $G_Y(v)$.
With this modification, the entire discussion in Example
3.2 adapts to this setting,
yielding 
\begin{equation} \label{poly-5}
C = \frac{1}{Ê\xi \lambda^\prime(\xi) {\bf E} \left[ \tau \right] }  {\bf E}_\xi \left[ \left( V_0 + \sum_{n=1}^\infty \frac{G_{Y_n}(V_{n-1})}{A_1\cdots A_n}
  \right)^\xi {\bf 1}_{\{ \tau = \infty\}}\right].
\end{equation}
As before $V_0 \sim \nu$,
where $\nu$ is obtained from the minorization $({\cal M})$, while
 $\tau$ is the first regeneration time and $\xi$ is the solution to the equation ${\bf E} \big[ A^\xi \big] = 1$.
The expectation is computed in 
the $\xi$-shifted measure, defined in the same manner 
as in \eqref{meas-shift},
but with $(A,B_0,\ldots,B_{k-1})$ in place of $(A,B,D)$.  
While this new representation can no longer be viewed as backward 
iterations of the original SFPE, it nonetheless
provides an algorithm with which the constant may be computed, and where the successive terms in the series become increasingly insignificant
as $n \to \infty$.

To see that this series is convergent and
to study the rate of convergence, 
using a slight variant of \eqref{pf-lm5.4.4},
it follows that the expectation in \eqref{poly-5} is bounded above by
\begin{equation} \label{poly-6}
\left({\bf E} \big[ V_0^\xi \big] + \sum_{n=1}^\infty {\bf E} \left[ \left(G_{Y_n}(V_{n-1})\right)^{\xi^\prime} \right]^{1/{\xi^\prime}}  {\bf P} \big( \tau > n \big)^{1/\xi^{\prime\prime}} \right)^\xi, 
\end{equation}
where $(1/\xi^{\prime}) + (1/\xi^{\prime\prime})=(1/\xi)$.

To estimate ${\bf E} \left[ \left(G_{Y_n}(V_{n-1}) \right)^\xi \right]$ 
as $n \to \infty$, it is helpful to first observe that
$\{ V_n \}$ converges to a stationary distribution.
Indeed,
using a technique of Loynes (\shortciteN{CG91}, p.\ 161), 
it follows that for any positive $\mbox{constant }c$,
\[
G_Y(v) \le \left( \sum_{i=1}^{k-1} B_i \right) c^{(k-1)/k} \max \left( 1, \frac{v}{c} \right) + B_0,
\]
and thus $F_Y(v) := Av + G_Y(v)$ satisfies
\begin{equation} \label{poly-7}
F_Y(v) \le \left( Ac + c^{(k-1)/k} \sum_{i=1}^{k-1} B_i  \right) \max  \left( 1, \frac{v}{c} \right) + B_0
:= \check{A} \max \left(c,v \right) + \check{B},
\end{equation}
where $\check{B} = B_0$ and $\check{A} = \left( A +  c^{-1/k} \sum_{i=1}^{k-1} B_i  \right)$.  The upper bound on the right-hand side is 
Letac's Model E, which itself converges to a stationary distribution under
the conditions of \shortciteN{CG91}, Proposition 6.1, namely under the assumption that $\xi > 0$ and that ${\bf E} \left[ \log (1 \vee \check{B}) \right]
< \infty$.  This latter condition is subsumed in the first of the
following hypotheses:
\begin{align} 
&{\bf E} \big[ B_i^{\xi^\prime} \big] < \infty,
\quad i=0,\ldots,k-1, \quad \mbox{\rm for some}\:\: \xi^\prime > \xi.
 \tag{$H_4$}\\
&{\bf E} \Big[ \log \Big( A + \Big(\sum_{j=1}^{k-1} B_j\Big)/2 B_0^\prime
\Big) \Big] < 0, \tag{$H_5$}
\end{align}
where $B_0^\prime$ is an independent copy of $B_0$.

Thus we conclude that the Letac model given on the right-hand side
of \eqref{poly-7} converges to its stationary distribution.
Next observe by $(H_5)$ and the argument on p.\ 162 of \shortciteN{CG91}
that the backward iterates of the 
SFPE \eqref{poly-1} converge and are
 independent of the initial value, and hence by Letac's principle (Lemma 2.1), the forward recursive sequence
also converges a.s.\ to the same random variable, which we denote by $V$.  Then $G_{Y_n}(V_n)$ converges a.s.\ to
\begin{equation} \label{poly-10}
G_Y(V) := \sum_{j=0}^{k-1} B_j V^{j/k},
\end{equation}
where $\{ B_j \}$ is independent of $V$ on the right-hand side.  

Finally, to obtain finiteness of the $\xi^\prime$th moment of \eqref{poly-10}, note that the process $\{ V_n \}$ is dominated
from above by
\begin{equation} \label{poly-15}
\check{V}_n := \check{A}_n \check{V}_{n-1} + \left(\check{A}_n c + \check{B}_n \right), \quad n=1,2,\ldots, \quad \check{V}_0 = V_0,
\end{equation}
which upon iteration yields
 $\check{V}_n = V_0 (\check{A}_1 \cdots \check{A}_n) \sum_{i=1}^n (\check{A}_i c + \check{B}_i )(\check{A}_{i+1} \cdots \check{A}_n)$.
But given $\alpha < \xi$, we can choose a constant $c$ such that ${\bf E} \big[ \check{A}^\alpha \big] <1$.
Then a simple computation (using either Minkowskii's inequality in the case
$\xi^\prime \ge 1$ or  \eqref{ex1-13} otherwise) yields that ${\bf E} \big[ \check{V}_n^\alpha \big]$ is uniformly bounded in $n$ for any $\alpha < \xi$.
Since $\{ \check{V}_n \}$ dominates $\{ V_n \}$, it now follows, using 
in \eqref{poly-10} the independence of $V$ and $\{B_0,\ldots,B_{k-1} \}$, that
\begin{equation} \label{poly-16}
\lim_{n \to \infty} {\bf E} \left[ \left( G_{Y_n}(V_{n-1})\right)^{\xi^\prime} \right] = \sum_{j=0}^{k-1} {\bf E} \big[ B_j^{\xi^\prime} \big] {\bf E} \left[ V^{\xi^\prime(j/k)} \right] < \infty
\end{equation}
provided that $\xi^\prime$ has been chosen sufficiently small
such that $\xi^\prime(k-1)/k < \xi$.  Then using the geometric recurrence
of $\{ V_n \}$, it follows that \eqref{poly-6} is finite; thus $C$
is finite.

For the above computations to make sense, it is critical that we can also verify the regularity properties of the Markov chain $\{ V_n \}$.
But $\varphi$-irreducibility and the minorization and drift conditions all follow by very minor modifications of the proofs given for Letac's Model E.
Furthermore, transience follows in the $\xi$-shifted measure, since under our assumptions $V_1 \ge B_0 > 0$,
and then $V_n \ge B_0(A_2 \cdots A_n) \uparrow \infty$ as $n \to \infty$.  
Thus we have established the following:

\begin{thm}  
Assume the polynomial model \eqref{poly-1}, and suppose that $(H_1)$,  $(H_4)$, and $(H_5)$ are satisfied.  Then
\beq{\label{eq2.1a}}
\lim_{u \ra \ff} u^{\xi}{\bf P} \left( V>u \right) =C,
\eeq
where $C$ is a finite positive constant given by \eqref{poly-5}.
\end{thm}

It is also straightforward to obtain an upper bound along
the lines of Theorem 2.2.
The methods outlined above can be modified to analyze a variety of recursions of the general form $V \stackrel{\cal D}{=} AV + G_Y(V)$
for an asymptotically  ``dominant" term $AV$ and ``remainder" term
$G_Y(V)$, where $v^{-1} G_Y(v) \downarrow 0$ as $v \to \infty$.  For
example, the composition $f \circ h$, where $f$ corresponds
to the SFPE in \eqref{poly-1} and $h$ is as in Letac's model \eqref{sr200},
can be handled in this manner.

%
%
%  MARKOV AND CONCLUSIONS
%
%
\section{Markov-driven processes and other extensions}
\subsection{Markov-driven processes and the ruin problem
with invesments}
Finally, we consider the case where the forward or backward
iterations in Letac's Model E are governed by an underlying Markov chain.
In this case, it is possible to first apply Lemma 2.2 to the {\em underlying
Markov chain}
to divide its excursions into independent blocks, and to observe
that the resulting process itself
has the form of Letac's Model E.  Thus, the theory developed in this paper may still be applied.

Several examples of this type have been considered in \shortciteN{JC09}, including 
Markovian generalizations of the ruin problem with investments, perpetuities,
and the GARCH(1,1) financial time series model. 

For example, in the ruin problem
with investments, the goal is to determine $\Psi(u) = 
{\bf P} \left( \sup_{n\in \pintegers} {\cal L}_n > u \right),$
where ({\it cf.} \eqref{intro3a})
\begin{equation} \label{Markovext1}
{\cal L}_n := B_1 + A_1 B_2 + (A_1 \cdots A_{n-1})B_n,
\end{equation}
and in a Markovian setting it is natural to assume, instead, that
$(A_n, B_n)= h(X_n)$, where $h$ is a deterministic function and $\{ X_n \}$
is a Harris recurrent Markov chain.  Now by Lemma 2.2, the chain $\{ X_n \}$
may be subdivided into independent blocks of the form 
$\{ X_{K_{i-1}}, \ldots, X_{K_i-1} \}$.  Consequently, upon setting
$J^\ast = J^\ast(n) := \inf\{i:  K_i >n \} -1$, we obtain
\begin{equation} \label{Markovext2}
{\cal L}_n = \check{B_0} + \check{A_0} \check{B}_1 
+ \cdots + (\check{A_0} \cdots \check{A}_{J^\ast-1})\check{B}_{J^\ast} 
+ (\check{A_0} \cdots \check{A}_{J^\ast}) {\cal R}_n,
\end{equation}
where
\begin{align*}
\check{A}_i& = A_{K_{i-1}} \cdots A_{K_i-1},\\ 
\check{B}_i & = B_{K_{i-1}} + A_{K_{i-1}} B_{K_{i-1}+1} + \cdots 
  + (A_{K_{i-1}} \cdots A_{K_i-2})B_{K_i-1},\\ 
\mbox{\rm and}\quad {\cal R}_n & = B_{K_{J^\ast}} + A_{K_{J^\ast}} B_{K_{J^\ast}+1} + \cdots
     + (A_{K_{J^\ast}} \cdots A_{n-1}) B_n
\end{align*}
for all $i=0,1,\ldots$, where $K_{-1}=1$.
Now by a short argument, given in \shortciteN{JC09}, it can be shown that if we
start the chain $\{ X_n \}$ at a regeneration time, then it follows from
\eqref{Markovext2} that ${\cal L} := \sup_{n \in \pintegers}{\cal L}_n$
satisfies an SFPE, namely,
\begin{equation} \label{Markovext3}
{\cal L} \stackrel{\cal D}{=} \max \left(\check{M},\check{A}{\cal L}\right)
 + \check{B},
\end{equation}
where $(\check{A},\check{B},\check{M}) \stackrel{\cal D}{=}
(\check{A}_i,\check{B}_i,\check{M}_i)$ and
\begin{align*}
\check{M}_i := \sup\{& B_{K_{i-1}} + A_{K_{i-1}}B_{K_{i-1}+1} + \cdots
 \\ 
 & + (A_{K_{i-1}} \cdots A_{j-1})B_j:  K_{i-1} \le j < K_i\} - \check{B}_i,
\quad \mbox{\rm for all}\:\: i.
\end{align*}
Finally, it can be verified
that---even without starting the process at a regeneration time---the 
asymptotics of ${\bf P} \left( {\cal L} >u \right)$ obtained via
the SFPE \eqref{Markovext3} are nonetheless the correct. 
For details, see \shortciteN{JC09}, Lemma 4.2.

The main conclusion to be drawn from the above discussion is that 
the Markov-driven process also satisfies an SFPE having the form of
Letac's Model E.  However, if the Markov chain takes values in a general
state space, as opposed to a finite state space, then it is not immediately
clear whether the required moment conditions in $(H_1)$ and $(H_2)$ 
of this paper will
actually be satisfied.  For the ruin problem with investments,
this issue has been addressed in detail in \shortciteN{JC09},
Theorem 4.2 and Section 6.1.  
Namely set $S_n = \sum_{i=1}^n \log A_i$,
and consider the G\"{a}rtner-Ellis limit from large deviation theory,
\[
\Lambda(\alpha) := \limsup_{n \to \infty} \frac{1}{n} \log
{\bf E} \left[ e^{\alpha S_n} \right], \quad \forall \alpha,
\]
and assume that the equation $\Lambda(\xi)=0$ has a unique positive solution.
Next, set $\tilde{S}_n = \sum_{i=1}^n \log A_i + \log |B_n|$ and
assume that both $\Lambda(\alpha)$ and 
\[
\tilde{\Lambda}(\alpha ):=  \limsup_{n \to \infty} \frac{1}{n} \log
{\bf E} \left[ e^{\alpha \tilde{S}_n} \right]
\]
are finite in a neighborhood of $\xi$.  Moreover, assume
that $\{ X_n \}$ satisfies the weak
regularity condition $(H_3)$ and minorization
condition $({\mathfrak M})$ given on pp.\ 1411-2 of \shortciteN{JC09}
(taking $h(x)=e^{\epsilon x}$ in ${(\mathfrak M})$, for any $\epsilon > 0$).
Together, these form one set of possible conditions 
suitable to apply the main results of \shortciteN{JC09}.
Then applying our current results to identify the constants of decay, we arrive
at the following extension of Theorem 2.1 of \shortciteN{JC09}.

\begin{thm}
Assume the conditions of \shortciteN{JC09}, and assume that $(\check{A},
\check{B},\check{M})$ satisfies the condition $(H_3)$ of this paper.
Then 
\begin{equation}
\Psi(u) \sim C u^{-\xi} \quad \mbox{\rm as}\:\: u \to \infty,
\end{equation}
where $C$ is given by \eqref{ex1-11}
but with
$\{ (A_i,B_i)  \}_{i \in \pintegers}$ replaced with
$\{ (\check{A}_i,\check{B}_i) \}_{i \in \pintegers}$, defined following
\eqref{Markovext2}, and where $\check{B}_0$
is independent with $\check{B}_0 \sim {\mathfrak L } (\check{B}^+)$.
\end{thm}

\subsection{Markov-driven Letac models}
The previous discussion can be extended by observing that
\eqref{Markovext1}
is the backward recursive sequence of the SFPE $f(v) = Av + B$,
and we could equally well study the supremum of a backward
sequence corresponding to a general Letac model, {\em viz.}
$f(v) = A \max\left(D,v\right) + B$.
As in the previous example, the limiting
perpetuity sequence for this more general model can be subdivided
into i.i.d.\ blocks.  Namely, letting $Z_n$ denote the $n$th
backward iterate of the process, and assuming that the driving
sequence for these backward iterates is of the form
$(A_n,B_n,D_n) = h(X_n)$ for a Harris recurrent Markov chain $\{ X_n \}$,
then just prior to regeneration, we have
\begin{equation} \label{markov3}
V_{K_i-1} = \max \left( \check{D}_i, \check{A}_i V_{K_{i-1}-1} \right) + \check{B}_i,
\end{equation}
where
\begin{align*}
\check{A}_i& = A_{K_{i-1}} \cdots A_{K_i-1}, \quad \check{B}_i =
\sum_{j=K_{i-1}}^{K_i-1} B_j \prod_{k=K_{i-1}}^{j-1} A_k,\\
& \hspace*{1cm}\mbox{\rm and} \quad \check{D}_i = \bigvee_{k=K_{i-1}}^{K_i-1} \left\{ D_k
  \prod_{j=K_{i-1}}^k A_j
 - \sum_{j=k+1}^{K_i-1} B_j \prod_{l=K_{i-1}}^{j-1} A_l \right\}.
\end{align*}
Setting $\check{Z} = \lim_{i \to \infty}Z_{K_i}$ 
and taking the limit as $i \to \infty$ yields the SFPE
\begin{equation}
\check{Z} \stackrel{\cal D}{=} \max \left( \check{D}, \check{A} \check{Z} \right)
+ \check{B},
\end{equation}
and under conditions analogous to those of Theorem 8.1, one naturally expects
this limiting tail distribution to correspond with that of $\lim_{n \to \infty} Z_n$.
We note that a slightly different SFPE will be obtained if, instead, one
considers the supremum of the process $\{ Z_n \}$ (involving a slight
modification of the random variable $\check{M}$ appearing in
\eqref{Markovext3}).  However, in either case, the resulting SFPE
has the same form as Letac's Model E, and hence the theory developed
in this paper applies.

Forward iteration can be handled similarly, although the random
variables $(\check{A},\check{B},\check{D})$ will be different
than those obtained from backward iteration.  Namely one obtains
\begin{equation} \label{markov3a}
V_{K_i-1} = \max \left( \check{D}_i, \check{A}_i V_{K_{i-1}-1} \right) + \check{B}_i,
\end{equation}
where
\begin{align*}
\check{A}_i& = A_{K_{i-1}} \cdots A_{K_i-1}, \quad \check{B}_i =
\sum_{j=K_{i-1}}^{K_i-1} B_j \prod_{k=j+1}^{K_i-1} A_k,\\
& \hspace*{1cm}\check{D}_i = \bigvee_{k=K_{i-1}}^{K_i-1} \left\{ D_k
  \prod_{j=k}^{K_i-1}A_j
 - \sum_{j=K_{i-1}}^{k-1} B_j \prod_{l=j+1}^{K_i-1} A_l \right\}.
\end{align*}
As a special case, forward iteration corresponds to a GARCH(1,1) process
with Markov modulated innovations; see Theorem 2.2 of \shortciteN{JC09},
where the corresponding moment and regularity conditions are studied
in a general Markovian setting for such processes.  Again, Theorem 2.1 of this
paper may be applied, although in all of these cases, it would be of interest to
determine whether the pair $(\check{A},\check{B})$ or triplet $(\check{A},\check{B},
\check{D})$ appearing in the expression for the constant $C$ can itself
be characterized probabilistically in a more transparent way.

\section{Concluding remarks}
In this paper we have introduced techniques that aid in the characterization of 
the constants $C$ and $\xi$  for solutions to SFPEs in the stationary case.  
Using the results of this paper, we have concurrently developed statistical methods for 
joint asymptotic inference  concerning $C$ and $\xi$.  This allows
for the development of inferential methods for such quantities
as Value at Risk and Expected Shortfall for GARCH (1,1) and related financial 
time series models, which are used extensively in applied econometrics. 
Extensions of our work to continuous time processes and general Markovian 
processes are presently under study. 

\bibliographystyle{chicago}
\bibliography{Bbl}

\end{document}